\documentclass[english]{article}
\usepackage{a4wide}
\usepackage{lmodern}
\usepackage[T1]{fontenc}
\usepackage[utf8]{inputenc}
\usepackage{geometry}
\usepackage{bm}
\usepackage{amsmath}
\usepackage{amssymb}
\usepackage{graphicx}
\usepackage{subfig}

\makeatletter

\newcommand{\lyxaddress}[1]{
\par {\raggedright #1
\vspace{1.4em}
\noindent\par}
}

\makeatother

\begin{document}

\renewcommand\bf[1]{\mathbf{#1}}
\newcommand*{\tr}{\mathsf{T}}
\renewcommand\[{\begin{equation}}
\renewcommand\]{\end{equation}}

\title{TraceFEM for the membrane problem using distance functions on $P_{1}$
and $P_{2}$ tetrahedra}

\author{Mirza Cenanovic}
\maketitle

\lyxaddress{\begin{center}
\textit{Department of Mechanical Engineering,  J\"onk\"oping University,
SE-55111  J\"onk\"oping, Sweden}
\par\end{center}}

\begin{abstract}
We consider Trace finite element methods for the linear membrane problem
on second order tetrahedral elements. To accomplish this, zero-level
set reconstruction methods for second order tetrahedra are considered.
For the higher order membrane model a corresponding stabilization
is proposed and numerically evaluated. We compare combinations of
background- and surface element order and provide numerical convergence
results. The impact of the stabilization on the resulting solution
is numerically analyzed. We also compare the choice of level set function
with respect to the geometrical distance and normal errors. 
\end{abstract}

\section{Introduction}

In this paper we extend the construction of finite element methods
for linear elastic membranes embedded in three dimensional mesh in
\cite{Cenanovic201698} to second order tetrahedral elements. 

We use the tangential calculus approach suggested for modeling surface
stresses in \cite{Gurtin1975}, for shells in \cite{Delfour1995}
and for finite element methods in \cite{Dziuk1988}. This approach
has recently become widely used see, e.g., \cite{Dziuk2013} for an
extensive overview. It was previously used on a triangulated surface
membrane in \cite{Hansbo2014}. 

We use a form of unfitted finite element approach suggested originally
in \cite{Olshanskii2009} where instead of using a triangulated representation
of the surface, the surface is implicitly defined on a background
mesh of higher dimension and the partial differential equations are
discretized on this mesh but integrated (or restricted) to the surface.
The surface is defined by the zero level-set of a level-set function.
This approach is also known as the TraceFEM and has become increasingly
popular recently, see e.g., \cite{Olshanskii2016} and the references
therein for an recent overview. One of the reasons why it is so attractive
from a numerical point of view lies in the way it handles moving (time
dependent) surfaces without the need for re-meshing techniques. Another
nice property of this approach is that complex shapes can be modeled
by implicit surfaces and directly used in simulations without the
need for costly mesh processing where often additional human interaction
is needed to clean up the computational mesh. Since the surface is
allowed to intersect the bulk mesh arbitrarily, small cuts will severely
affect the resulting conditioning of the linear systems. Thus, we
adapt a ghost penalty stabilization approach proposed in \cite{BuHaLa15}
and used for a variety of different surface and bulk-surface problems,
e.g., \cite{Hansbo2015,Cenanovic201698,Cenanovic2015,Burman2016,Hansbo2017}.
In this paper, we adapt the ghost penalty approach for a second order
TraceFEM for the membrane problem. Development of stabilization methods
for TraceFEM is currently a hot topic and other stabilization methods
exist, such as the full gradient stabilization. See e.g., \cite{Burman2016a}. 

In previous works in \cite{Cenanovic2015,Cenanovic201698} we used
the name CutFEM for this method. The name CutFEM however is more suited
for methods where the bulk solution is used, for pure surface problems
we prefer the name TraceFEM as it becomes more clear what is implied.

Recent focus has been put into the development of methods for the
reconstruction and numerical integration of implicitly defined domains,
see e.g., \cite{Olshanskii2016} for a recent overview. In this work
we adapt the approach suggested by \cite{Fries2015,Fries2017} where
the idea is to interpolate the implicit function using a standard
parametric interpolation of order $m$ and employ a Newton-Raphson
root finding algorithm to reconstruct the zero-level set geometry. 

\subsection{Overview}

This work is divided as follows. We begin by introducing the membrane
model and its TraceFEM in Section \ref{sec:MembraneModel}. The details
of reconstructing a second order implicit surface are explained in
Section \ref{sec:Zero-Level-Surface-reconstruction}. The resulting
numerical error estimations are presented in Section \ref{sec:Numerical-Results}.
Finally Section \ref{sec:Conclusion} provides a conclusion and discussion
about future work.

\section{Membrane model and Finite Element Method\label{sec:MembraneModel}}

\subsection{Tangential calculus}

Let $\Gamma$ denote a smooth surface which is embedded in $\mathbb{R}^{3}$
and has an outward pointing normal $\bm{n}_{\Gamma}$. The surface
contains two types of boundaries, $\partial\Gamma_{\mathrm{N}},$
where we assume zero traction boundary conditions, and $\partial\Gamma_{\text{D}},$
were we assume zero Dirichlet boundary conditions. 

\begin{figure}
\centering{}\includegraphics[width=0.75\textwidth]{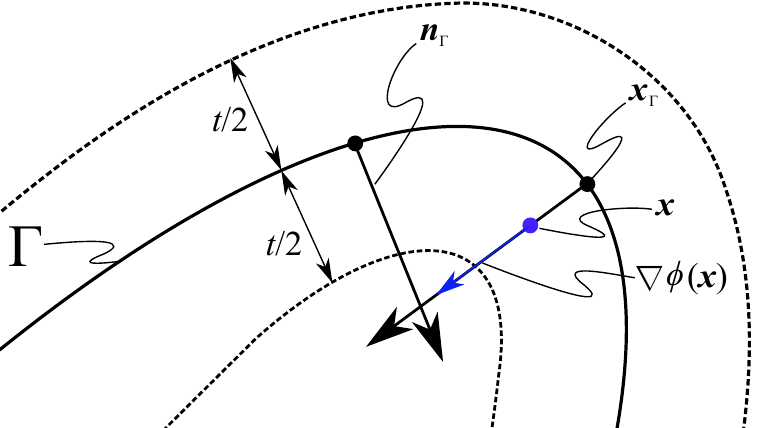}\caption{Tubular neighborhood of $\Gamma$ defined by the dashed line. The
surface is said to be sufficiently smooth when the normal vectors
$\nabla\phi(\bm{x}_{\Gamma})$ do not intersect within the tubular
neighborhood, in other words each point $\bm{x}$ in this neighborhood
has an unique closest point $\bm{p}(\bm{x})$ to $\Gamma$.\label{fig:Tubular-neighborhood}}
\end{figure}

We let $\phi(\bm{x})$ denote the signed distance to $\Gamma$ at
each point $\bm{x}\in\mathbb{R}^{3}$ and note that on $\bm{x}\in\Gamma$
the normal coincides with the gradient of the distance function $\nabla\phi|_{\bm{x}_{\Gamma}}=\bm{n}(\bm{x}_{\Gamma})$.
The domain that is occupied by the membrane is defined by

\[
\Omega_{t}=\{\bm{x}\in\mathbb{R}^{3}:|\phi(\bm{x})|<t/2\},
\]
where $t$ is the thickness of the membrane, see Figure \ref{fig:Tubular-neighborhood}.
Using a signed distance function we ensure that $\text{|\ensuremath{\nabla\phi}|=1}$
and then, given a sufficiently smooth surface, we can assume that
a function $u$ on $\Gamma$ can be extended to the neighborhood of
$\Gamma$ by means of a closest point projection $\bm{p}(\bm{x})=\bm{x}-\phi(\bm{x})\nabla\phi(\bm{x})$
such that $u^{e}(\bm{x})=u(\bm{p}(\bm{x}))$ . 

 The surface gradient $\nabla_{\Gamma}$ on $\Gamma$ is defined by

\begin{equation}
\nabla_{\Gamma}u=\bm{P}_{\Gamma}\nabla u^{e},\label{eq:TangentialGradient}
\end{equation}
where $\nabla$ denotes the full $\mathbb{R}^{3}$ gradient and $\bm{P}_{\Gamma}=\bm{P}_{\Gamma}(\bm{x})$
the orthogonal projection of $\mathbb{R}^{3}$ onto the tangential
plane of $\Gamma$ at $\bm{x}\in\Gamma$ given by

\begin{equation}
\bm{P}_{\Gamma}=\bm{I}-\bm{n}\otimes\bm{n},\label{eq:TangentialProjection}
\end{equation}
where $\bm{I}$ is the identity matrix. It is readily shown that the
tangential gradient (\ref{eq:TangentialGradient}) is independent
of the extension $u^{e}$ (see, e.g., \cite[Chapter 9, Section 5]{Delfour2011b}),
hence no distinction will be made between functions on $\Gamma$ and
their extensions in what follows.

The surface gradient and its components are denoted by

\[
\nabla_{\Gamma}u=:\left(\frac{\partial u}{\partial x_{\Gamma}},\frac{\partial u}{\partial y_{\Gamma}},\frac{\partial u}{\partial z_{\Gamma}}\right).
\]

The tangential Jacobian matrix for a vector valued function $\bm{v}(\bm{x})$
is defined as the dyadic product of $\nabla_{\Gamma}$ and $\bm{v}$,

\renewcommand*{\arraystretch}{2}

\[
(\nabla_{\Gamma}\otimes\bm{v})^{\tr}:=\begin{bmatrix}\dfrac{\partial v_{1}}{\partial x_{\Gamma}} & \dfrac{\partial v_{1}}{\partial y_{\Gamma}} & \dfrac{\partial v_{1}}{\partial z_{\Gamma}}\\
\dfrac{\partial v_{2}}{\partial x_{\Gamma}} & \dfrac{\partial v_{2}}{\partial y_{\Gamma}} & \dfrac{\partial v_{2}}{\partial z_{\Gamma}}\\
\dfrac{\partial v_{3}}{\partial x_{\Gamma}} & \dfrac{\partial v_{3}}{\partial y_{\Gamma}} & \dfrac{\partial v_{3}}{\partial z_{\Gamma}}
\end{bmatrix},
\]

\renewcommand*{\arraystretch}{1}

\noindent and the surface divergence is $\nabla_{\Gamma}\cdot\bm{v}:=\mathrm{tr}(\nabla_{\Gamma}\otimes\bm{v})$.
For the vector valued function $\bm{u}$ the surface strain tensor
is defined by

\[
\bm{\varepsilon}(\bm{u}):=\dfrac{1}{2}\left(\nabla_{\Gamma}\otimes\bm{u}+(\nabla_{\Gamma}\otimes\bm{u})^{\tr}\right)
\]

\noindent and the in-plane strain tensor is defined by

\[
\bm{\varepsilon}_{\Gamma}(\bm{u}):=\bm{P}_{\Gamma}\bm{\varepsilon}(\bm{u})\bm{P}_{\Gamma}.
\]

\subsection{The membrane model}

Following \cite{Hansbo2014,Cenanovic201698}, we consider the problem
of finding $\bm{u}:\Gamma\rightarrow\mathbb{R}^{3}$ such that

\begin{equation}
\begin{aligned}-\nabla_{\Gamma}\cdot\bm{\sigma}_{\Gamma}(\bm{u})=\bm{f} & \quad\text{on }\Gamma,\\
\bm{\sigma}_{\Gamma}=2\mu\bm{\varepsilon}_{\Gamma}+\lambda\text{tr}(\bm{\varepsilon}_{\Gamma})\bm{P}_{\Gamma} & \quad\text{on }\Gamma,\\
\bm{u}=0 & \quad\text{on }\partial\Gamma_{D},\\
\bm{\sigma}\cdot\bm{n}=0 & \quad\text{on }\partial\Gamma_{N},
\end{aligned}
\label{eq:membraneProblem}
\end{equation}

\noindent where $\bm{f}:\Gamma\rightarrow\mathbb{R}^{3}$ is an area
load, 

\[
\mu:=\dfrac{E}{2(1+\nu)},\quad\lambda:=\dfrac{E\nu}{1-\nu^{2}}
\]

\noindent are the Lamé parameters in plane stress where $E$ denotes
the Young's modulus and $\nu$ Poisson's ratio. Under the assumption
that the material obeys Hooke's law under plane stress, these equations
can be derived from the minimization of the surface potential energy
equation

\[
\Pi_{\Gamma}(\bm{u}):=\dfrac{1}{2}\int_{\Gamma}\bm{\sigma}_{\Gamma}(\bm{u}):\bm{\varepsilon}_{\Gamma}(\bm{u})d\Gamma-\int_{\Gamma}\bm{f}\cdot\bm{u}d\Gamma
\]

\noindent as shown in \cite{Cenanovic201698}. The weak form of (\ref{eq:membraneProblem})
is defined by: find $\bm{u}\in V$ such that

\[
a(\bm{u},\bm{v})=l(\bm{v}),\ \forall\bm{v}\in V,
\]

\noindent where

\begin{flalign*}
a(\bm{u},\bm{v}) & =\left(2\mu\bm{\varepsilon}_{\Gamma}(\bm{u}),\bm{\varepsilon}_{\Gamma}(\bm{v})\right)_{\Gamma}+\left(\lambda\nabla_{\Gamma}\cdot\bm{u},\nabla_{\Gamma}\cdot\bm{v}\right)_{\Gamma}\\
 & =\left(2\mu\bm{\varepsilon}(\bm{u}),\bm{\varepsilon}(\bm{v})\right)_{\Gamma}-\left(4\mu\bm{\varepsilon}(\bm{u})\cdot\bm{n},\bm{\varepsilon}(\bm{v})\cdot\bm{n}\right)_{\Gamma}+\left(\lambda\nabla_{\Gamma}\cdot\bm{u},\nabla_{\Gamma}\cdot\bm{v}\right)_{\Gamma},
\end{flalign*}

\[
l(\bm{v})=(\bm{f},\bm{v})_{\Gamma},
\]

\noindent and

\[
(\bm{v},\bm{w})_{\Gamma}=\int_{\Gamma}\bm{v}\cdot\bm{w}d\Gamma,\quad\text{and}\quad\left(\bm{\varepsilon}_{\Gamma}(\bm{v}),\bm{\varepsilon}_{\Gamma}(\bm{w})\right)_{\Gamma}=\int_{\Gamma}\bm{\varepsilon}_{\Gamma}(\bm{v}):\bm{\varepsilon}_{\Gamma}(\bm{w})d\Gamma
\]

\noindent are the $L_{2}$ inner products.

\subsection{The trace finite element method}

This section describes the discretization using TraceFEM. Let $\tilde{\mathcal{K}}_{h}$
denote a quasi uniform mesh into shape regular tetrahedra of a domain
$\Omega$ in $\mathbb{R}^{3}$ completely containing $\Gamma$. In
this work we define the surface $\Gamma$ implicitly by constructing
a signed continuous scalar distance function $\phi(\bm{x})$ such
that

\[
\Gamma=\{\bm{x}\in\Omega:\phi(\bm{x})=0\},
\]

\noindent which is a continuous zero-isosurface. It should be noted
that the property $|\nabla\phi|=1$ does not need to hold in general,
i.e., it is not necessary for $\phi$ to be a distance function in
the actual computations; however, if it holds, then the zero-isosurface
becomes less sensitive to small perturbations. It can also be beneficial
in cases of evolving surfaces, cf. \cite{Cenanovic2015}. 

\begin{figure}

\centering{}\subfloat[]{\begin{centering}
\includegraphics[width=0.45\textwidth]{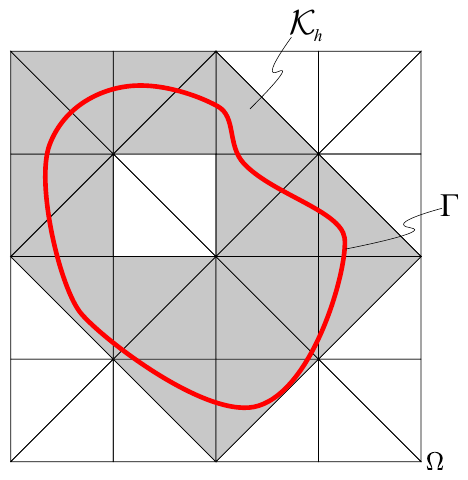}
\par\end{centering}
}\subfloat[]{\centering{}\includegraphics[width=0.4\textwidth]{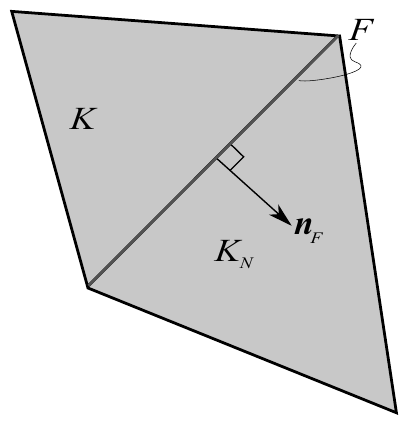}}\caption{a) Active background mesh $\mathcal{K}_{h}$. b) Interior face $F$
shared by two active background elements $K$ and $K_{N}$.\label{fig:BackgroundMesh}}
\end{figure}

The active background mesh is defined as the set of background elements
that are cut by the zero-isosurface by

\[
\mathcal{K}_{h}=\{K\in\tilde{\mathcal{K}}_{h}:\tilde{\mathcal{K}}_{h}\cap\Gamma\neq\emptyset\}
\]
and its set of interior faces by

\[
\mathcal{F}_{h}=\{F=K\cap K_{N}:K,K_{N}\in\mathcal{K}_{h}\}.
\]
For all active cut elements $K\in\mathcal{K}_{h}$ there is a neighbor
$K_{N}\in\mathcal{K}_{h}$ such that $K$ and $K_{N}$ share a face,
see Figure \ref{fig:BackgroundMesh}. 

Let $\partial\Omega_{h,D}$ denote the boundary of discrete domain
$\Omega_{h}$ that is intersected by the discrete surface boundary
denoted as $\partial\Gamma_{h,D}$. The finite element space is then
defined by

\[
V_{h}=\left\{ \bm{v}\in[\tilde{V}_{h}|_{\Omega_{h}}]^{3}:\bm{v}=\bm{0}\text{ on }\partial\Omega_{h,D}\right\} 
\]
where $\tilde{V}_{h}$ is a space of continuous polynomials of order
$m_{B}=\{1,2\}$ (subscript $B$ denotes the bulk) defined on $\tilde{\mathcal{K}}_{h}$.
In this work, zero boundary conditions are treated by assuming that
$\partial\Omega_{h,D}$ intersects $\partial\Omega_{h}$, which is
accomplished simply by prescribing the displacements in the nodes
of the background mesh. For a more general handling of boundary conditions
we could use Nitsche's method, see, e.g., \cite{Burman2015,Lehrenfeld2016}.

The finite element method on $\Gamma_{h}$ is given by: find $\bm{u}_{h}\in V_{h}$
such that 

\[
A_{h}(\bm{u}_{h},\bm{v})=l_{h}(\bm{v})\ \forall\bm{v}\in V_{h},
\]
where the bilinear form $A_{h}(\cdot,\cdot)$ is defined by

\begin{gather*}
A_{h}(\bm{u}_{h},\bm{v}):=a_{h}(\bm{v},\bm{w})+\\
\begin{cases}
\gamma j_{h,1}(\bm{v},\bm{w}) & \text{if }m_{B}=1\\
\gamma_{1}j_{h,1}(\bm{v},\bm{w})+\gamma_{2}j_{h,2}(\bm{v},\bm{w}) & \text{if }m_{B}=2
\end{cases}\\
\forall\bm{v},\bm{w}\in V_{h}
\end{gather*}
with 
\[
a_{h}(\bm{v},\bm{w})=\left(2\mu\bm{\varepsilon}_{\Gamma_{h}}(\bm{v}),\bm{\varepsilon}_{\Gamma_{h}}(\bm{w})\right)_{\Gamma_{h}}+\left(\lambda\nabla_{\Gamma_{h}}\cdot\bm{v},\nabla_{\Gamma_{h}}\cdot\bm{w}\right)_{\Gamma_{h}},
\]

\[
j_{h,1}(\bm{v},\bm{w})=\sum_{F\in\mathcal{F}_{h}}\left(\left[\nabla\bm{v}\right],\left[\nabla\bm{w}\right]\right)_{F},
\]
and

\[
j_{h,2}(\bm{v},\bm{w})=\sum_{F\in\mathcal{F}_{h}}h^{2}\left(\left[\nabla^{2}\bm{v}\right],\left[\nabla^{2}\bm{w}\right]\right)_{F}.
\]
Here $j_{h}(\cdot,\cdot)$ denotes the face stabilization term, where
\[
[\nabla\bm{v}]=\nabla\bm{v}|_{K_{i}\bigcap F}-\nabla\bm{v}|_{K_{j}\bigcap F}
\]
and

\[
[\nabla^{2}\bm{v}]=\nabla^{2}\bm{v}|_{K_{i}\bigcap F}-\nabla^{2}\bm{v}|_{K_{j}\bigcap F}
\]
denotes the jump of $\nabla\bm{v}$ and $\nabla^{2}\bm{v}$ respectively
across $F$ and $\gamma$, $\gamma_{1}$ and $\gamma_{2}$ are scalar
stabilization parameters that are user defined. The discrete surface
gradients are defined using the normals to the discrete surface

\[
\nabla_{\Gamma_{h}}\bm{v}=\bm{P}_{\Gamma_{h}}\nabla\bm{v}=(\bm{I}-\bm{n}_{h}\otimes\bm{n}_{h})\nabla\bm{v},
\]
and the right hand side is given by

\[
l_{h}(\bm{v})=(\bm{f},\bm{v})_{\Gamma_{h}}.
\]

The face stabilization term $j_{h}(\cdot,\cdot)$ is used to reduce
ill-conditioning in $a_{h}(\cdot,\cdot)$, which results from the
surface arbitrary cutting through the background elements. The discrete
normals in case of $m>1$ are given by

\[
\bm{n}_{h}:=\dfrac{\dfrac{\partial\bm{x}_{\Gamma_{h}}}{\partial r}\times\dfrac{\partial\bm{x}_{\Gamma_{h}}}{\partial s}}{\left|\dfrac{\partial\bm{x}_{\Gamma_{h}}}{\partial r}\times\dfrac{\partial\bm{x}_{\Gamma_{h}}}{\partial s}\right|},
\]
where $\frac{\partial\bm{x}_{\Gamma_{h}}}{\partial r}$ is defined
using a parametric map $\bm{F}:(r,s)\rightarrow(x,y,z)$ on a reference
2D element $T^{r}$. Details on how to compute $\Gamma_{h}$ are given
in Section \ref{sec:Zero-Level-Surface-reconstruction}.

\section{Zero-level surface reconstruction\label{sec:Zero-Level-Surface-reconstruction}}

\begin{figure}
\centering{}\includegraphics[width=1\textwidth]{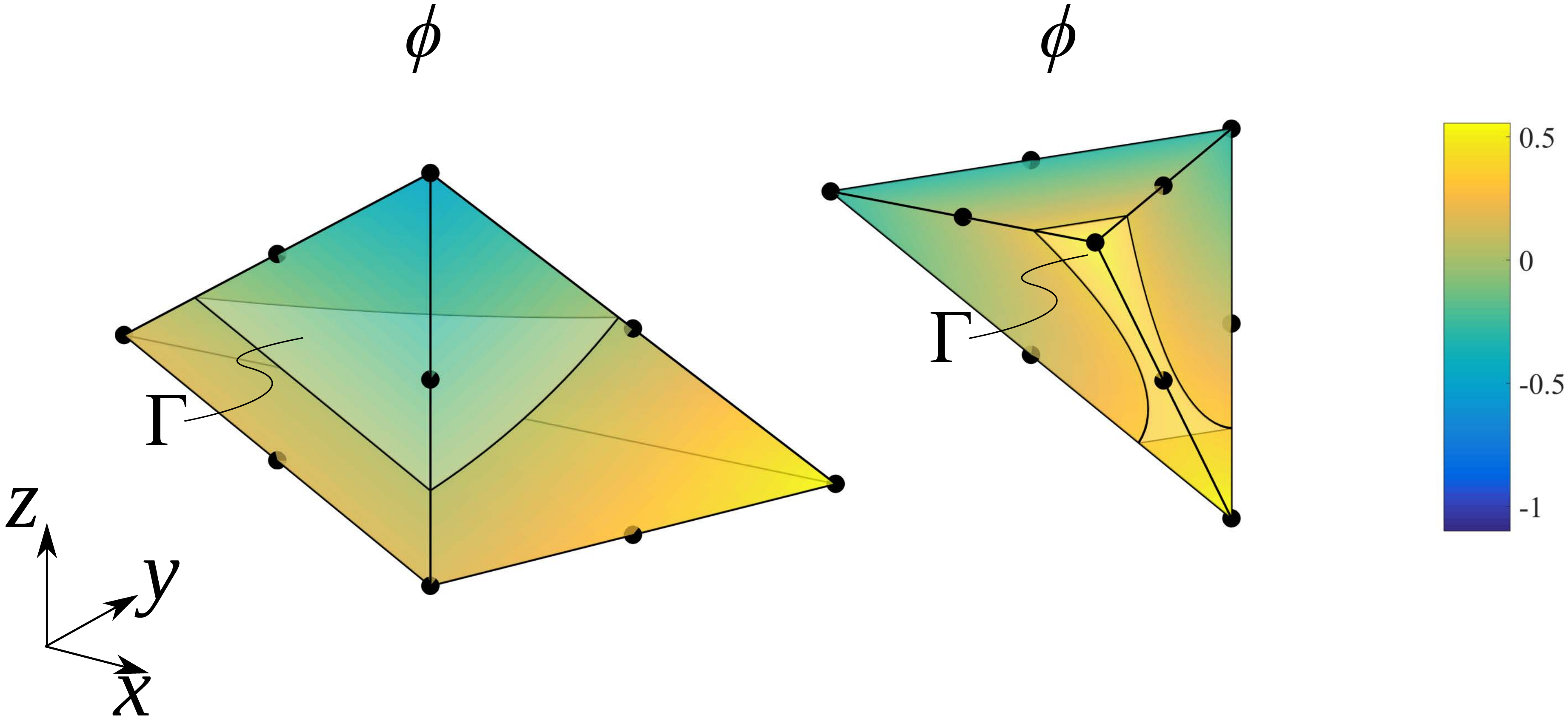}\caption{Isocontours of the continuous signed distance function $\phi$ on
a 10-noded tetrahedral element.\label{fig:Isocontours}}
\end{figure}

\begin{figure}
\centering{}\subfloat[]{\centering{}\includegraphics[width=0.45\textwidth]{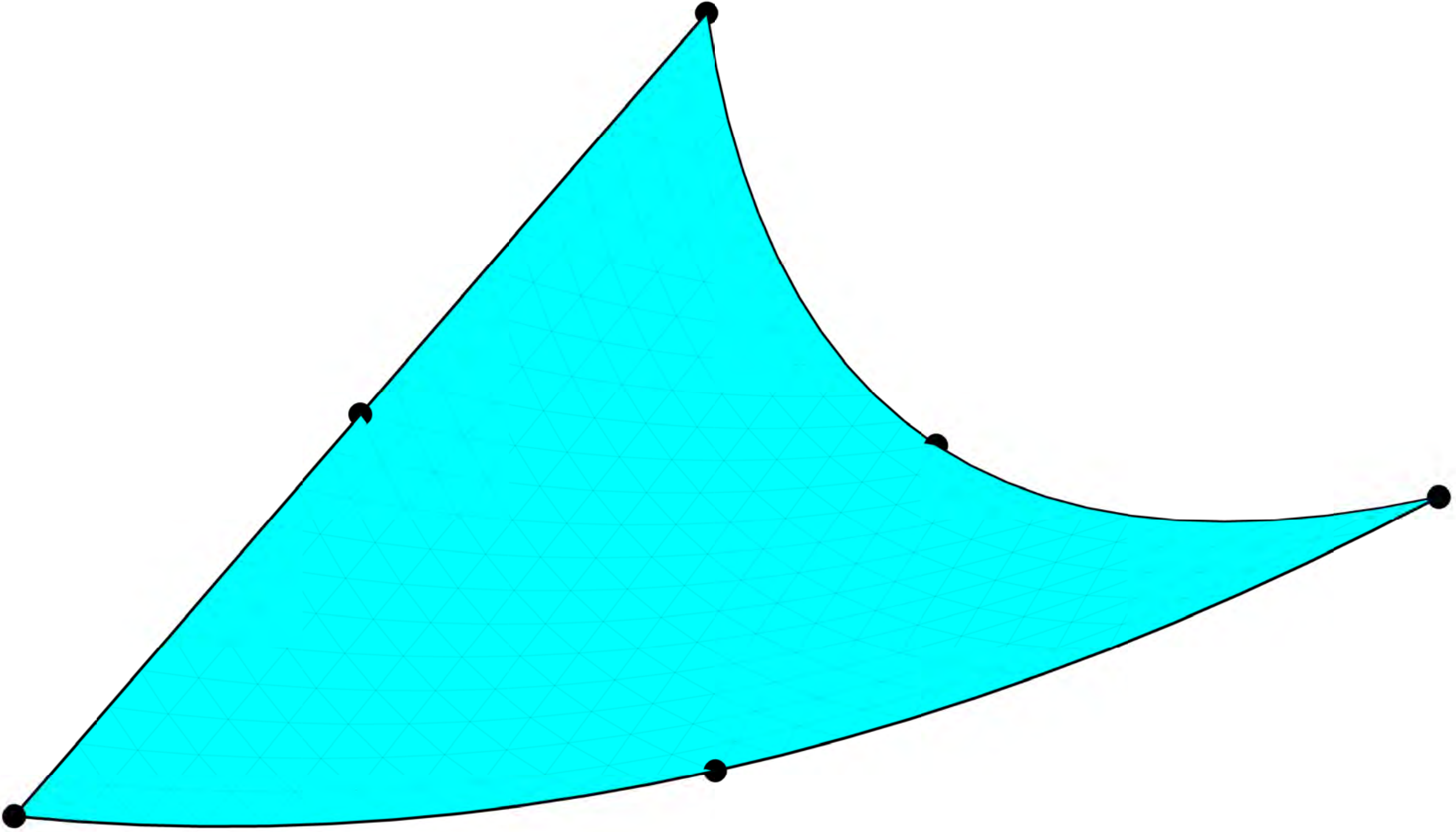}}$\quad$\subfloat[]{\centering{}\includegraphics[width=0.45\textwidth]{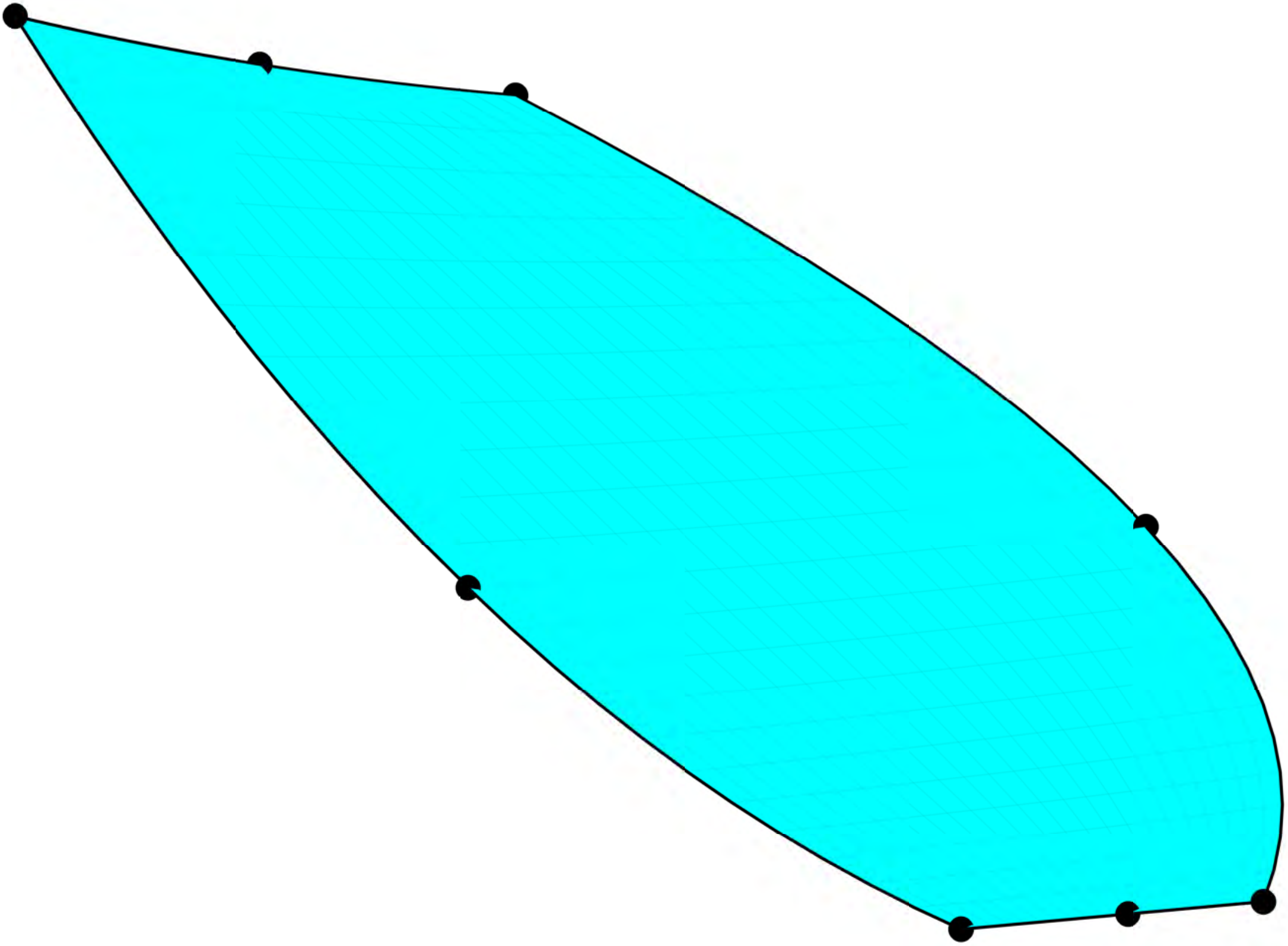}}\caption{Resulting surface types from cutting a tetrahedral element arbitrarily
with a zero-level set. (a) 6 extracted surface points mapped to a
$P_{2}$ triangle element. (b) 8 extracted surface points mapped to
a 8-noded serendipity element. \label{fig:Resulting-surface-types}}
\end{figure}

In this section we describe the approach for extracting the discrete
zero-level set $\Gamma_{h}$ from a signed distance function $\phi(\bm{x})$.
The basic idea is to determine the zero-level set for each element
$K$ by some form of root finding. In previous works \cite{Cenanovic2015,Cenanovic201698,Hansbo2015,Burman2015}
this was done by simple linear interpolation on linear tetrahedral
element. In these cases, the value of $\phi$ is exact in the nodes
of $K$ and the zero-level set $\Gamma_{h}$ is interpolated linearly
along the edges of $K$ yielding the corners of a planar surface element
$T$. Here, however, we need to find the zero-level points along the
edges of a second order tetrahedral element $K$ and the zero-level
points that lie on the faces of $K$, see Figure \ref{fig:Isocontours}.
The set of these surface points will define the nodes of a second
order surface Lagrange element. Since $K$ is a second order tetrahedral
element, which is assumed to be affine, the arbitrary intersection
with a surface will yield two types of surface elements; second order
triangles and second order quadrilaterals, see Figure \ref{fig:Resulting-surface-types}.
In general, a continuous $\phi(\bm{x})$ might not be known, instead
we may only have access to a discrete signed distance function $\bm{\Phi}$
defined in the nodes of $K$. In this case we can create an approximation
of $\phi$using of the basis functions of the bulk element $K$:

\[
\phi_{h}(\bm{x})=\sum_{i\in N_{K}}\varphi_{i}^{m_{B}}(\bm{x})\Phi_{i},
\]
where $N_{K}$ is the set of nodes in $K$, $\Phi_{i}$ are the known
nodal values of the signed distance function and $\varphi_{i}^{m_{B}}(\bm{x})$
is the basis function of polynomial order $m_{B}=\{1,2\}$ acting
on element $K$. Note that the basis functions can alternatively be
mapped or defined in the physical coordinate system since the bulk
element is assumed affine. 

\subsection{Root finding}

Following the recent work done in \cite{Fries2015,Fries2017} we set
up methods for extracting the zero-level set from both the continuous
level set function $\phi(\bm{x})$ (if available) and the discrete
$\bm{\Phi}$ (which can always be available in the nodes of the background
mesh). We begin by denoting the two discrete zero-level sets

\[
\Gamma_{h}|_{\phi}=\{\bm{x}\in\Omega:\Pi_{h}^{m_{\Gamma}}\phi(\bm{x})=0\},
\]

\[
\Gamma_{h}|_{\phi_{h}}=\{\bm{x}\in\Omega:\Pi_{h}^{m_{\Gamma}}\phi_{h}(\bm{x})=0\},
\]
where the interpolant $\Pi_{h}^{m_{\Gamma}}$ of order $m_{\Gamma}$
is described presently. 

\subsubsection{Valid topology\label{subsec:Valid-topology}}

\begin{figure}
\begin{centering}
\subfloat[]{\centering{}\includegraphics[width=0.49\textwidth]{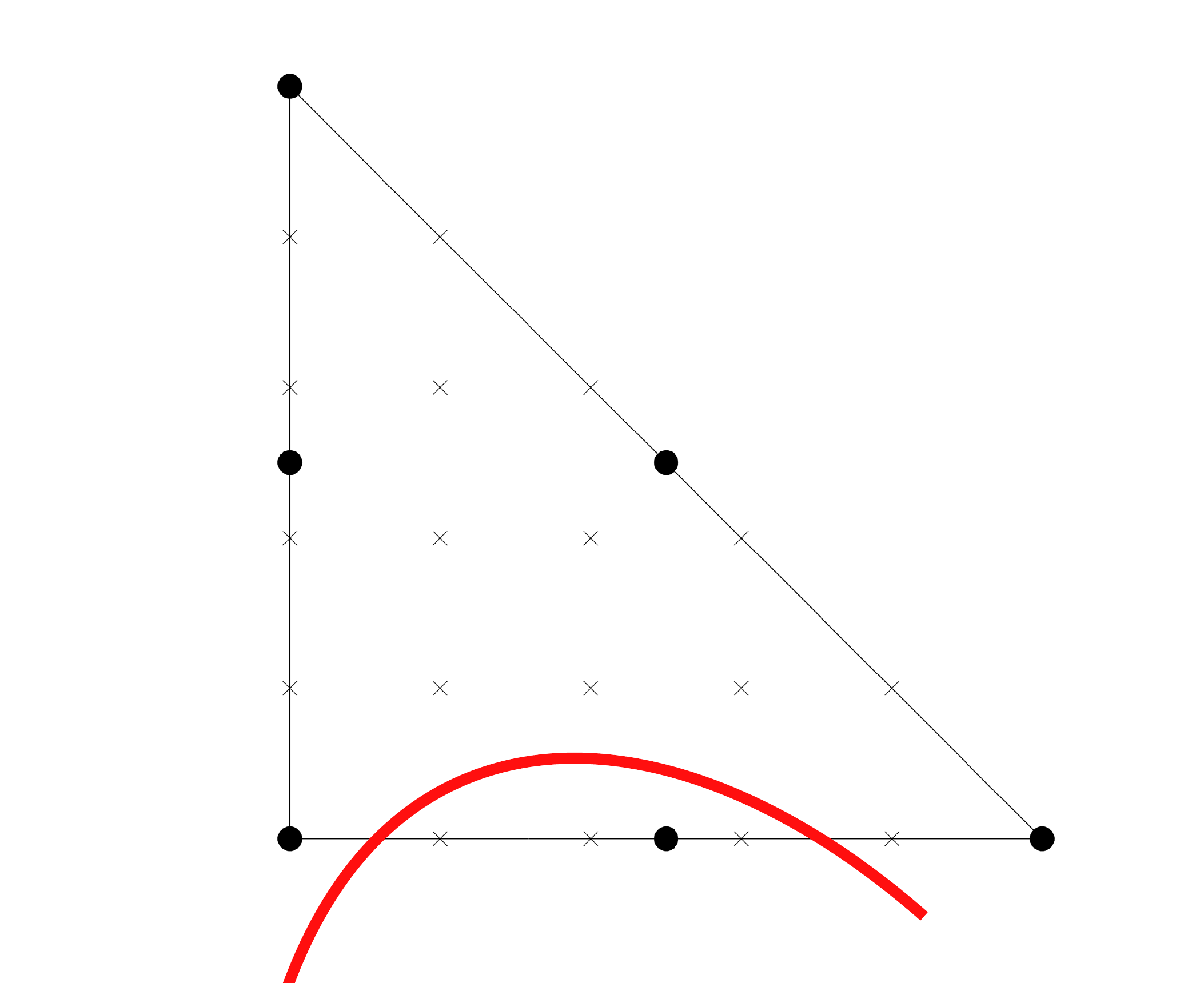}}\subfloat[]{\centering{}\includegraphics[width=0.49\textwidth]{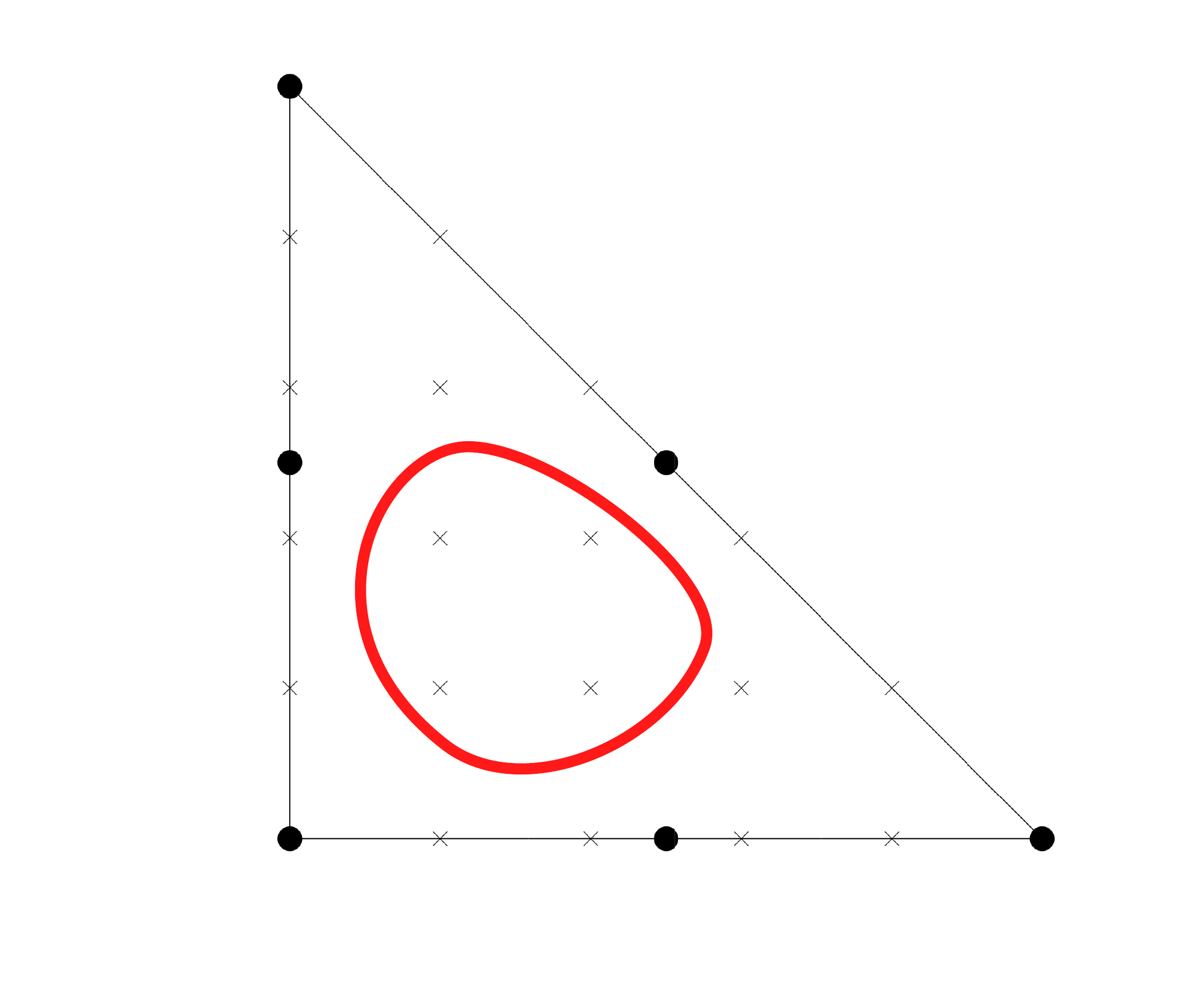}}
\par\end{centering}
\centering{}\subfloat[]{\centering{}\includegraphics[width=0.49\textwidth]{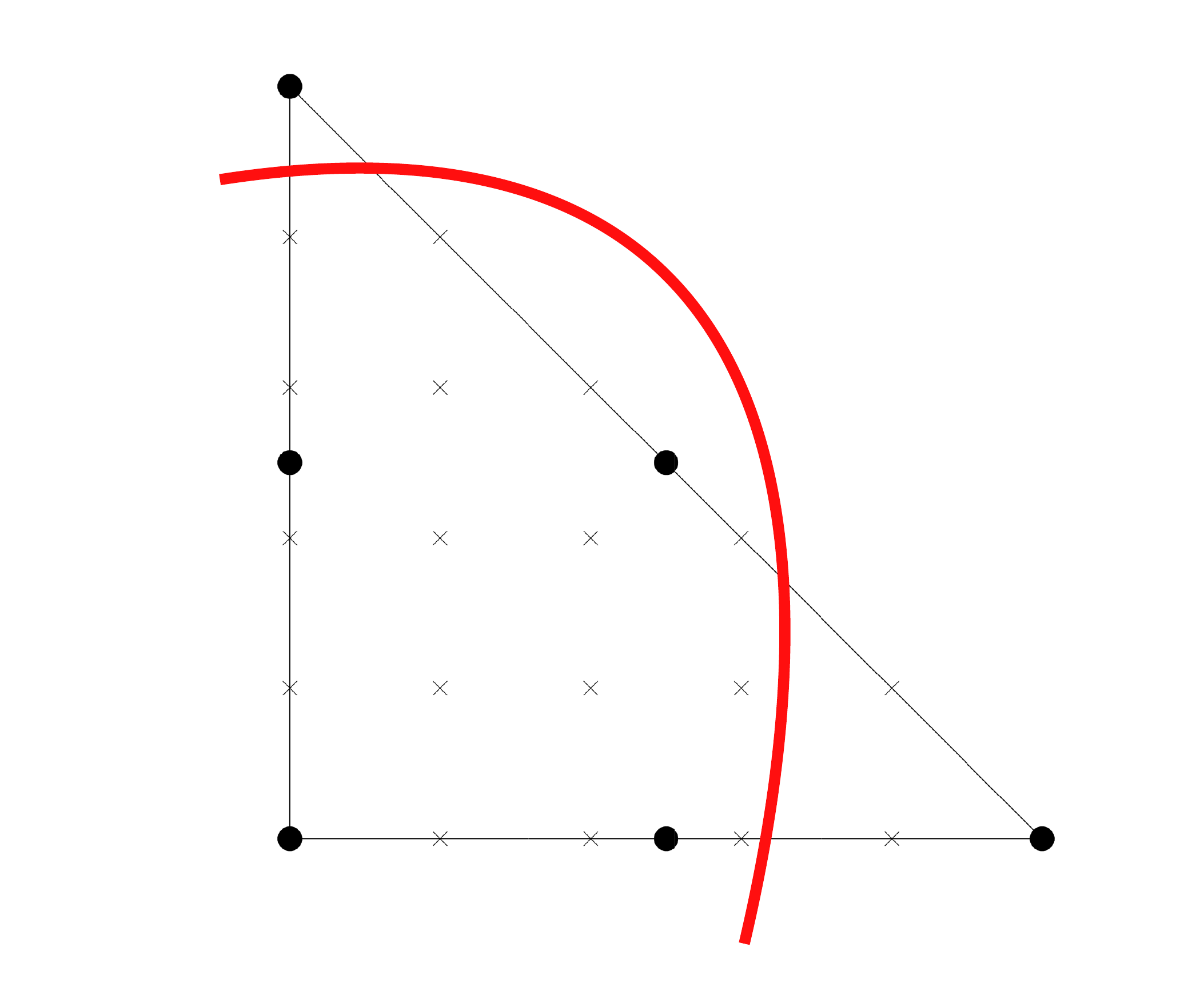}}\subfloat[\label{fig:BadTopoHighCurv}]{\centering{}\includegraphics[width=0.49\textwidth]{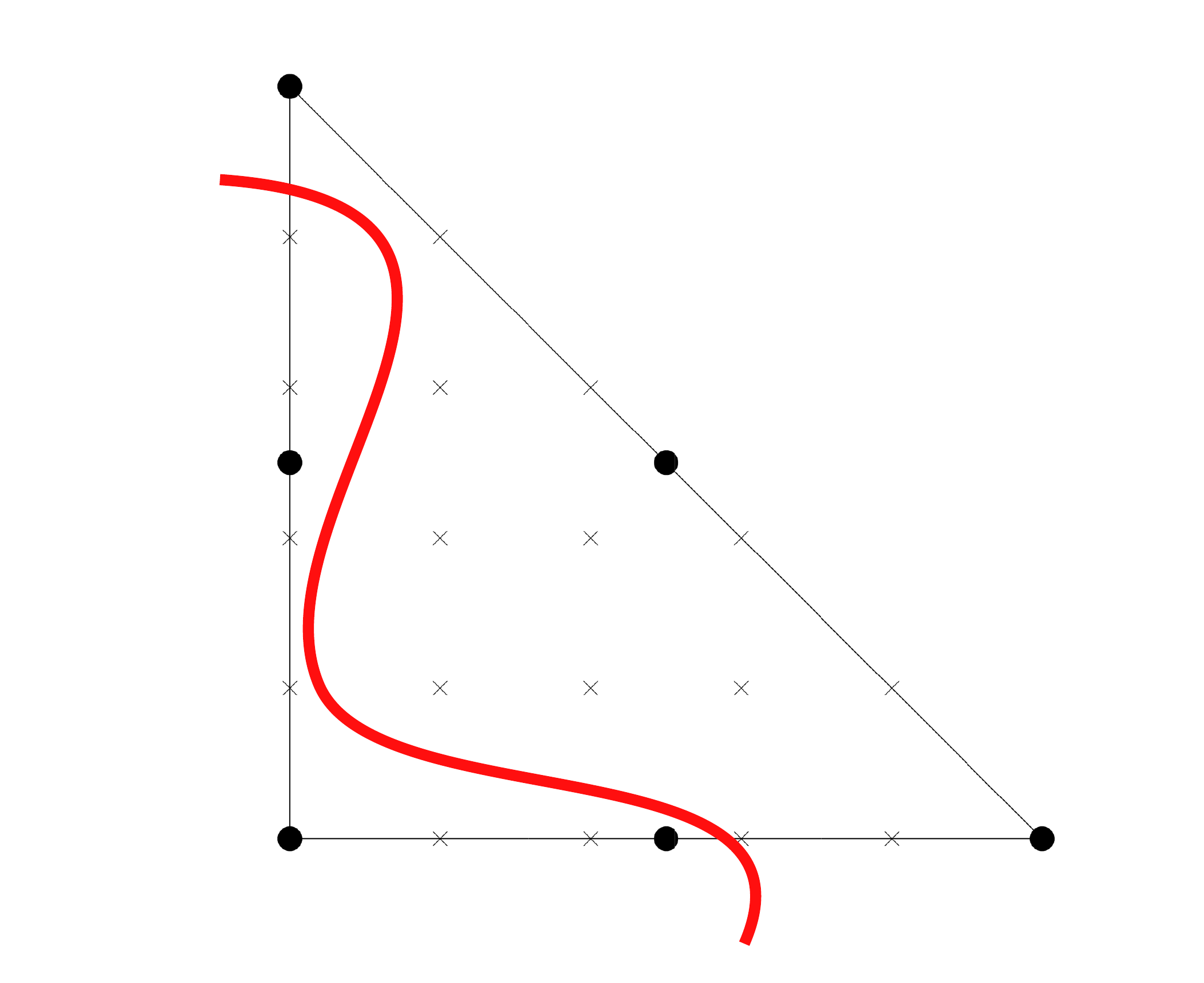}}\caption{Examples of bad topologies visualized on side-views of a 3D parametric
second order tetrahedral element with uniformly distributed sample
points. The red curves represent a surface. a) One edge cut more than
once. b) A small interface exists inside an element. c) More than
two edges are cut. d) The curvature of the cut is too high. \label{fig:BadTopology}}
\end{figure}

As is pointed out in \cite{Fries2015}, comparing the different signs
of the function $\bm{\Phi}|_{K}$ in the corner nodes of the element
$K$ is not sufficient to determine if the surface topology is valid.
Here a valid topology means that the arbitrary intersection of an
implicit surface with the faces of a background element result in
a number of surface points that can be mapped to polygons. To determine
if an element is cut we compute

\begin{equation}
\underset{i\in N_{\mathrm{grid}}}{\min}\left(\Phi_{i}^{\mathrm{grid}}\right)\cdot\underset{i\in N_{\mathrm{grid}}}{\max}\left(\Phi_{i}^{\mathrm{grid}}\right)<0,\label{eq:IsCut}
\end{equation}
where
\[
\Phi_{j}^{\mathrm{grid}}=\sum_{i=1}\varphi_{i}^{m_{B}}\left(\bm{r}_{j}^{\mathrm{grid}}\right)\cdot\Phi_{i}\quad\forall j\in N_{\mathrm{grid}},
\]
$\bm{r}_{j}^{\mathrm{grid}}$ denotes a number of uniformly spaced
sample points in the parametric space, see Figure \ref{fig:BadTopology}.
Note that $\varphi_{i}^{m_{B}}\left(\bm{r}_{j}^{\mathrm{grid}}\right)$
can be computed in a pre-processing step, and re-used for every background
element. For an example of bad topology see Figure \ref{fig:BadTopology}.
In the case described in Figure \ref{fig:BadTopoHighCurv} we identify
high curvature by 

\[
\nabla\bar{\Phi}^{\mathrm{grid}}\cdot\nabla\Phi_{j}^{\mathrm{grid}}<tol,
\]
where $\nabla\bar{\Phi}^{\mathrm{grid}}$ denotes the average of all
$\nabla\Phi_{j}^{\mathrm{grid}}$ for all $j\in N_{\mathrm{grid}}$
and $tol$ is a user defined number chosen such that large differences
in the angle between $\nabla\bar{\Phi}^{\mathrm{grid}}$ and $\nabla\Phi_{j}^{\mathrm{grid}}$
define high curvature, here $\nabla\Phi_{j}^{\mathrm{grid}}$ is given
by 

\[
\nabla\Phi_{j}^{\mathrm{grid}}=\sum_{i=1}\nabla\varphi_{i}^{m_{B}}\left(\bm{r}_{j}^{\mathrm{grid}}\right)\cdot\Phi_{i}.
\]

To be certain that the surface topology is valid we check the following
conditions on each face of the tetrahedral:
\begin{itemize}
\item Each edge of the face may only be cut once.
\item The number of cuts per face must be two.
\item If no face is cut, then all nodes of the tetrahedron must have the
same sign and thus the whole tetrahedron is uncut.
\end{itemize}
In case of invalid topology, local refinement can be used to resolve
the background mesh.

\subsubsection{The case of discrete level set function}

\begin{figure}
\centering{}\includegraphics[width=0.5\textwidth]{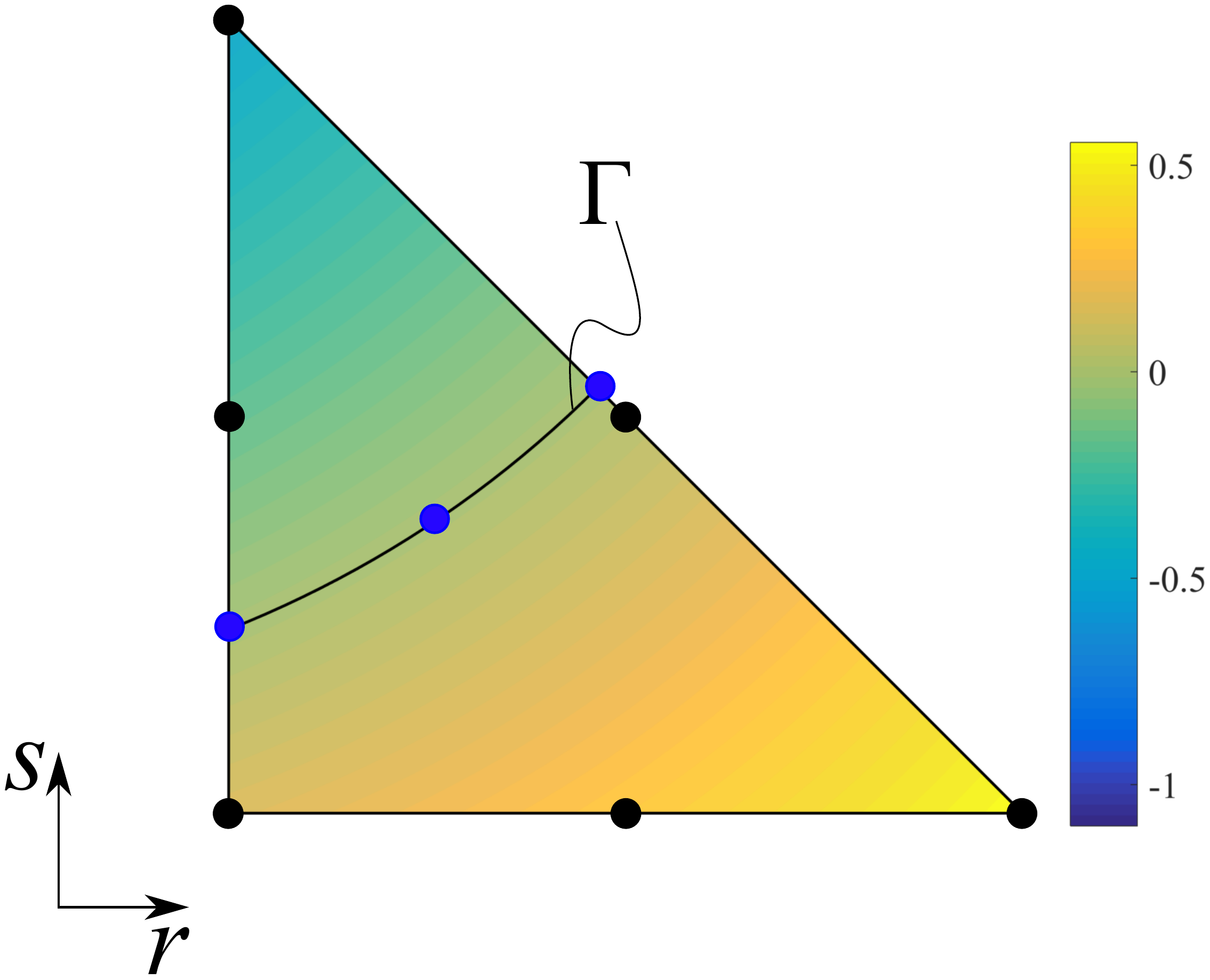}\caption{Isocontours of the distance function $\phi$ on a mapped face of a
tetrahedral element.\label{fig:IsocontoursParametric}}
\end{figure}

\begin{figure}
\centering{}\subfloat[]{\centering{}\includegraphics[width=0.45\textwidth]{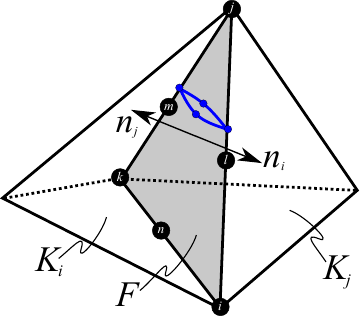}}\subfloat[]{\centering{}\includegraphics[width=0.45\textwidth]{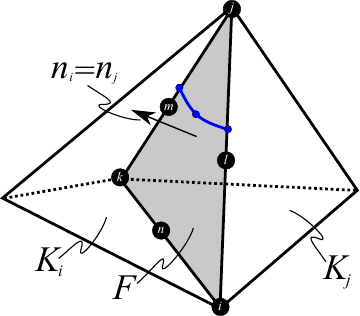}}\caption{Face numbering. (a) $F_{K_{i}}=\{i,j,k,l,m,n\}$, $F_{K_{j}}=\{i,k,j,n,m,l\}$.
(b) $F_{K_{i}}=F_{K_{j}}=\{i,j,k,l,m,n\}$.\label{fig:Face-numbering}}
\end{figure}

In order to find the roots for $\phi_{h}(\bm{x})=0$, when $\phi_{h}(\bm{x})$
is interpolated using an interpolant $\Pi_{h}^{m_{\Gamma}}$ of order
$m_{\Gamma}=2$, we follow the work done in \cite{Fries2015,Fries2017}
using the following steps:
\begin{enumerate}
\item For each element check if the topology is valid by following the steps
in the previous Section \ref{subsec:Valid-topology}. 
\item For each face $F$ of tetrahedral element $K$ in $\mathcal{K}_{h}|_{m_{B}=2}$
the nodal values of $\bm{\Phi}|_{F}$ are mapped to a parametric triangle
$T^{\text{r}}|_{m_{\Gamma}=2}$, see Figure \ref{fig:IsocontoursParametric}.
If element $K$ has a valid topology it's faces must have either two
or zero cut edges, additionally at least three faces must be cut.
We determine if the face is cut and which edges are cut by following
a procedure analogues to (\ref{eq:IsCut}). Additionally we renumber
the nodes of the faces such that they are unique, i.e., the normal
of each face $F_{K_{i}\cap K_{j}}$, no matter which tetrahedral they
belong to, points in the same direction, $\bm{n}|_{F_{K_{i}\cap K_{j}}}=\bm{n}|_{F_{K_{j}\cap K_{i}}}$.
This ensures that the gradients computed with the shape functions
of the face elements are the same for both elements $K_{i}$ and $K_{j}$,
otherwise the edge-points of the surface elements might not coincide,
see Figure \ref{fig:Face-numbering}. The resulting surface is thus
guarantied to be $C_{0}$ continuous.
\item On each cut edge on the parametric face $T_{m_{\Gamma}}^{r}$ we employ
a Newton-Raphson iterative search scheme:\\
\[
\bm{r}_{i+1}=\bm{r}_{i}-\dfrac{\phi_{h}(\bm{r}_{i})}{\nabla\phi_{h}(\bm{r}_{i})\cdot\bm{s}}\bm{s},
\]
where $\bm{r}=[r,s]$ is the local coordinate of the parametric triangle,
$\nabla\phi_{h}(\bm{r}_{i})$ is evaluated by interpolation using
the basis functions and $\bm{s}$ is the search direction. To find
the root along the edges, $\bm{s}$ is simply the directional vector
along the edge.
\item Once the two edge points are found the inner node needs to be determined
by the same root finding scheme. It turns out that the search direction
is critical for the convergence of the Newton search as well as the
geometrical convergence as shown in \cite{Fries2015,Fries2017}, where
the authors propose five different variations of the search directions
and two ways of starting position of the search. Choosing a linearly
interpolated starting position (straight line between the edge roots)
and set the search direction to be the normal to the line or $\bm{s}=\nabla\phi_{h}(\bm{r}_{0})$
yields satisfying results with respect to accuracy and performance,
see \cite{Fries2017}. In some rare cases when the Newton search fails
if gets stuck in a false root lying outside of the triangle, in this
case we employ bisection in order to get back inside the triangle
where the Newton search is continued until convergence. This approach
yields a robust method in all cases but increases the number of iterations
slightly for these rare cases.
\item The resulting surface points need to be numbered such that their normal
is oriented in the same general direction as $\nabla\bar{\Phi}^{\mathrm{grid}}$. 
\item Using this method we either get 6 surface points which are mapped
to a second order triangular element, or 8 points in which case we
map them to an 8-noded serendipity element. In the case of quadrilaterals
the reason for mapping to an 8-noded serendipity element is to avoid
the additional iterative search for the midpoint. Our argument against
splitting it into two triangles is that we get less integration points
which makes integration less expensive compared to two triangular
elements. 
\item If the resulting discrete surface needs to be used for smooth surface
shading, then an additional step is needed to create a connectivity
from the list of unconnected surface elements. In order to accomplish
this efficiently the background mesh information for each surface
patch is used to uniquely number the nodes and create the connectivity
list. Note that this step is not necessary for integration. 
\end{enumerate}
If we have access to the exact function $\phi$, the procedure above
is still valid, with the difference that we need to map $\bm{r}_{i}$
to $\bm{x}$ before evaluating $\phi(\bm{x}(\bm{r}_{i}))$ and $\nabla\phi(\bm{x}(\bm{r}_{i}))$. 

\begin{figure}
\centering{}\includegraphics[width=0.45\textwidth]{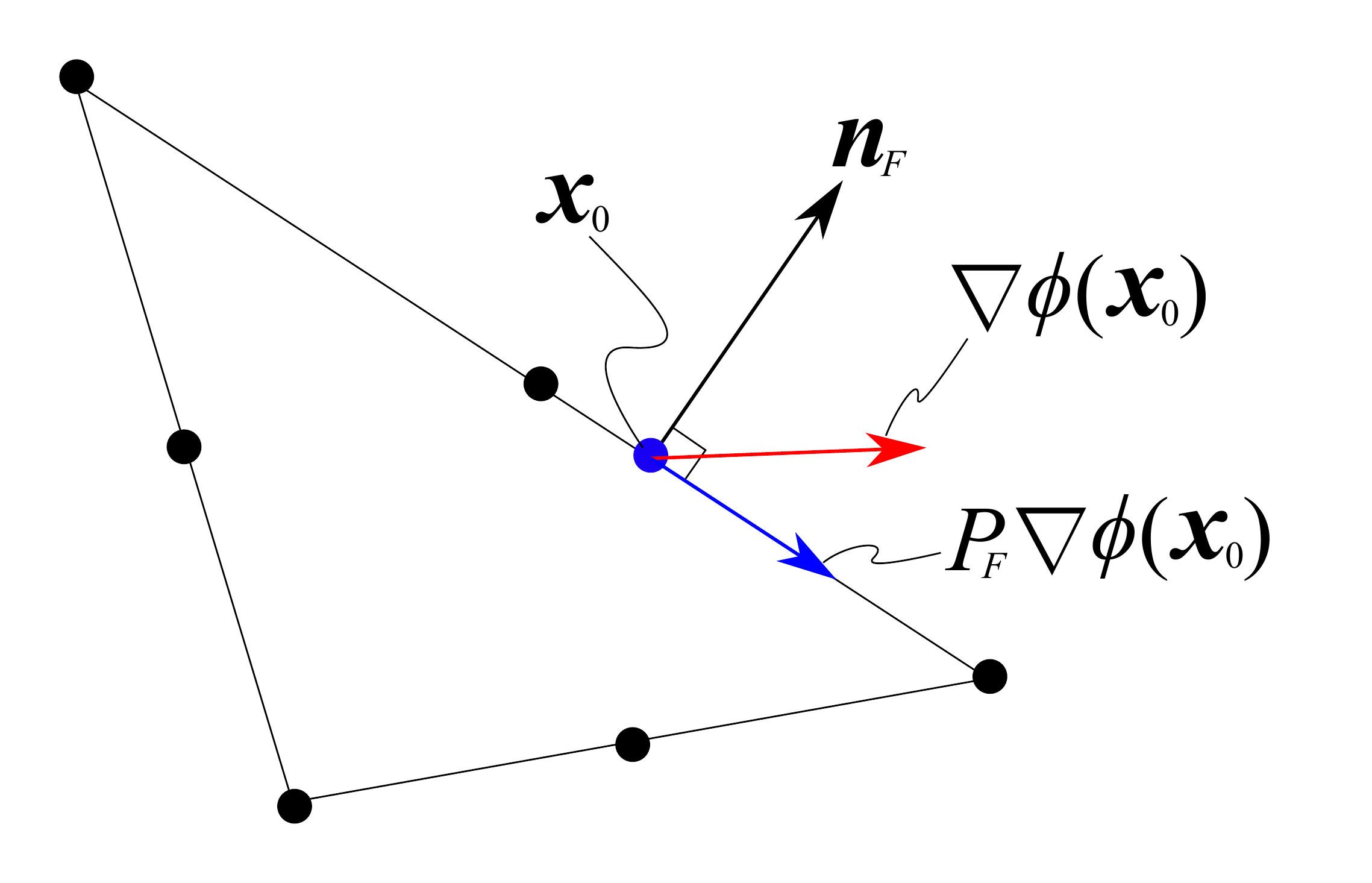}\caption{Sideview of a tetrahedral element. The search direction $\nabla\phi$
is projected onto the tetrahedral face $f$ (shown here as a line)
resulting in a modified Newton method with the search direction $P_{F}\nabla\phi$.
\label{fig:NewtonSearchAlongFace}}
\end{figure}

It is possible to create the above scheme in physical coordinates
by evaluating the basis functions in physical coordinates, $\varphi(\bm{x})$,
see the Appendix. The search for roots on edges in this case is the
same as above, the search on faces however is ``free'' since $\bm{s}=\nabla\phi(\bm{x}_{0})$.
In this case we restrict $\bm{s}$ to the (planar) face of the tetrahedron
by tangential projection:

\[
\bm{s}=\bm{P}_{F}\nabla\phi(\bm{x}_{0}),
\]
where $\bm{P}_{F}=\bm{I}-\bm{n}_{F}\otimes\bm{n}_{F}$, $\bm{I}$
is the identity matrix and $\bm{n}_{F}\otimes\bm{n}_{F}$ the outer
product of the face normal to the tetrahedron face, see Figure \ref{fig:NewtonSearchAlongFace}.
Note that the construction of $\bm{\varphi}(\bm{x})$ is done to avoid
the mapping of $\bm{r}$ to $\bm{x}$ in each step of the root finding
algorithm.

\section{Numerical Results\label{sec:Numerical-Results}}

The mesh size parameter for subsequent convergence studies is defined
as

\[
h:=\dfrac{1}{\sqrt[3]{\mathrm{N}}},
\]
where $\mathrm{N}$ is the number of nodes in the uniformly refined
mesh $\mathcal{K}_{h}.$ We denote the order of the surface elements
as $m_{\Gamma}$ and the order of the bulk elements as $m_{B}$. In
the tables the columns named ``Rate'' denote the rate of convergence.

\begin{figure}
\begin{centering}
\subfloat[$\Gamma_{h}|_{\phi}$]{\centering{}\includegraphics[width=0.49\textwidth]{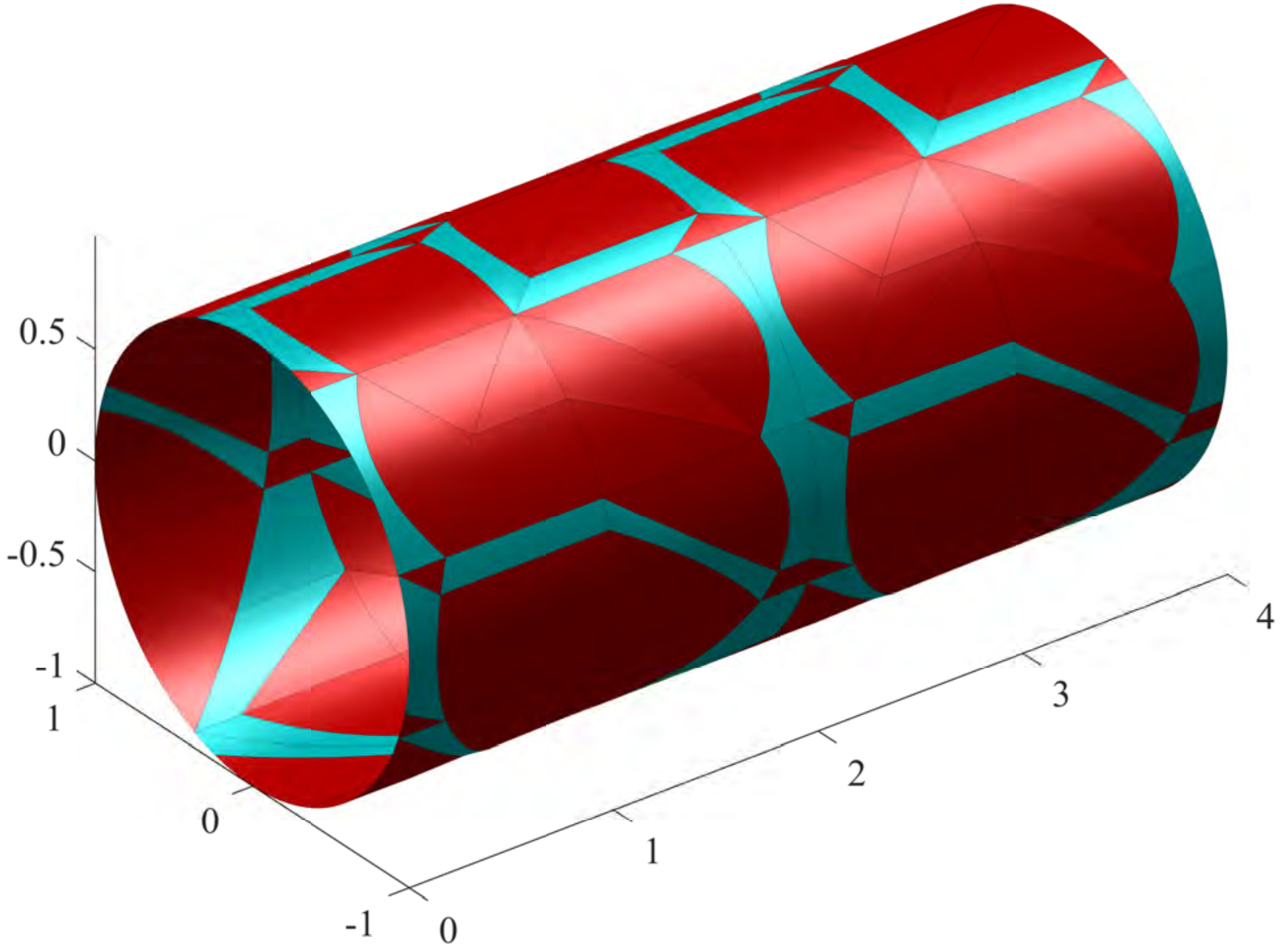}}\subfloat[$\Gamma_{h}|_{\phi_{h}}$]{\centering{}\includegraphics[width=0.49\textwidth]{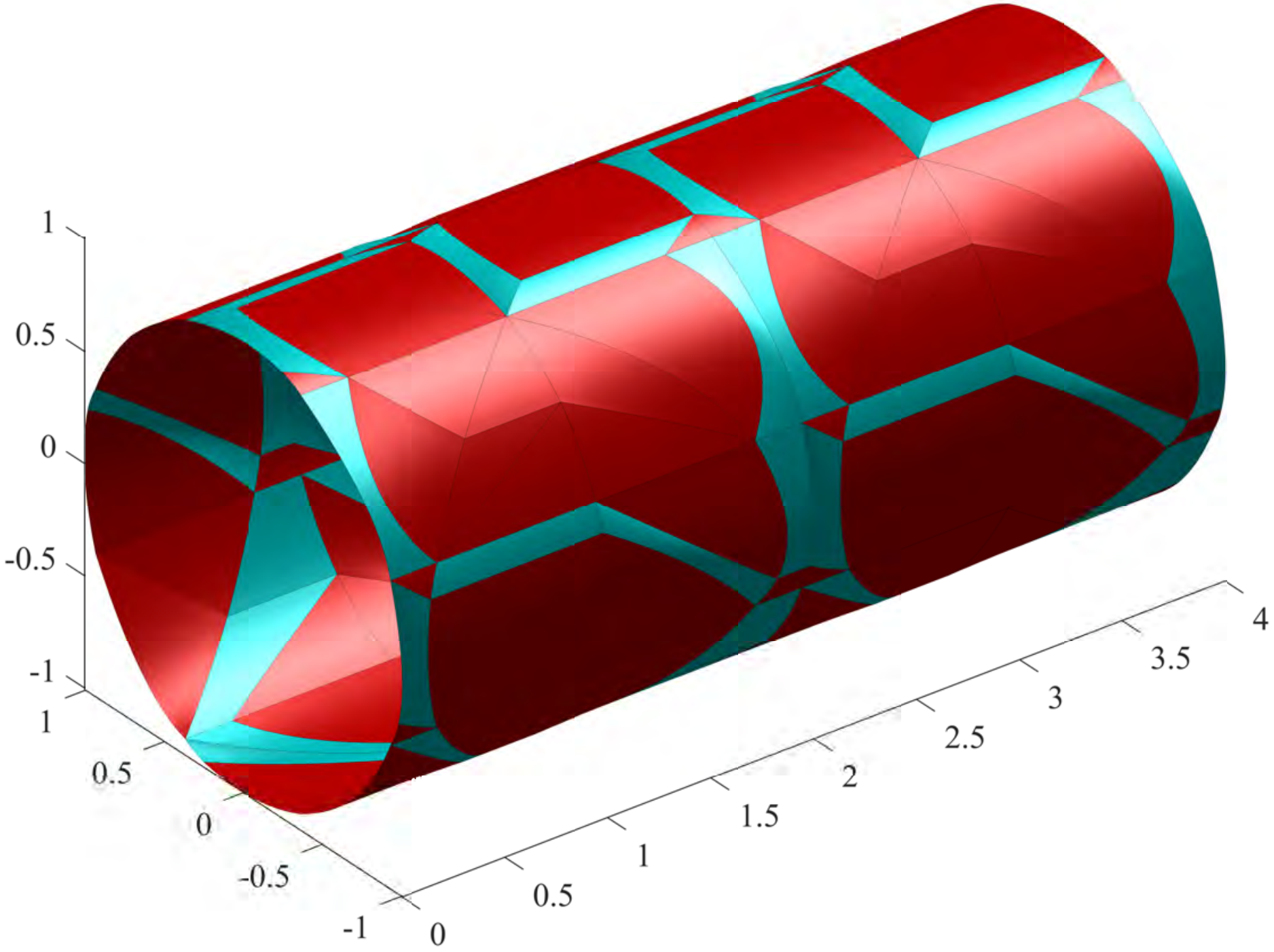}}
\par\end{centering}
\centering{}\subfloat[Front view. $\Gamma_{h}|_{\phi}$]{\centering{}\includegraphics[width=0.49\textwidth]{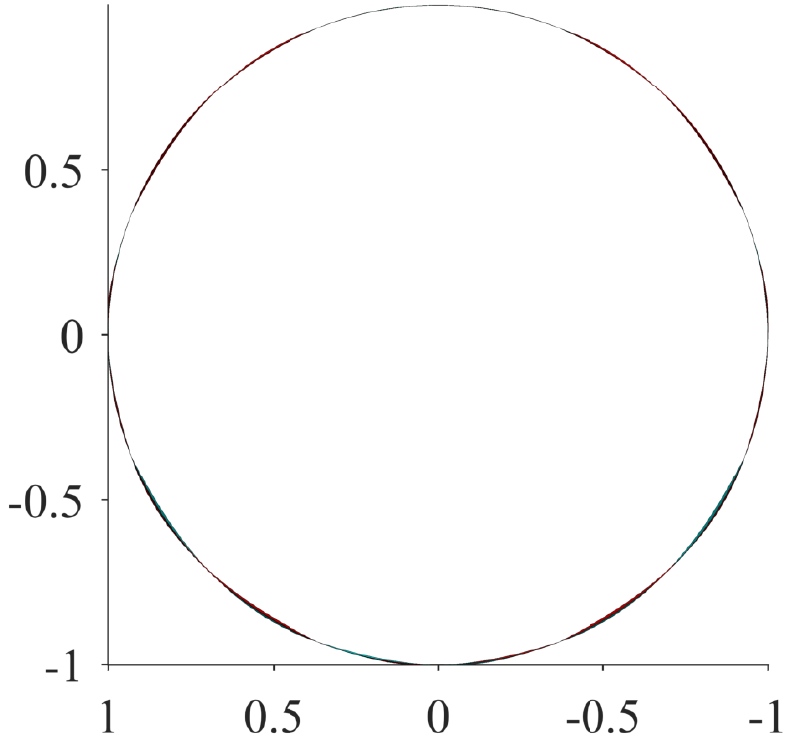}}\subfloat[Front view. $\Gamma_{h}|_{\phi_{h}}$]{\centering{}\includegraphics[width=0.49\textwidth]{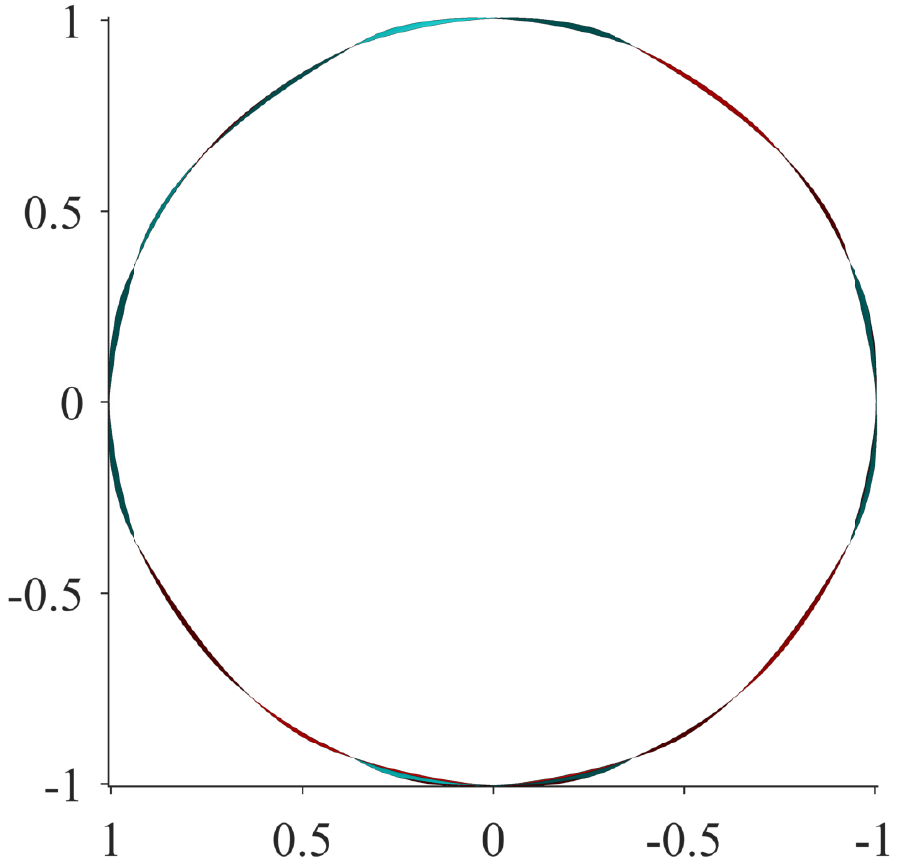}}\caption{Reconstructed surface with $m_{B}=2$.\label{fig:Reconstructed-surface}}
\end{figure}

The resulting reconstructed surfaces can be seen in Figure \ref{fig:Reconstructed-surface}.

\subsection{Membrane error comparison}

\begin{figure}
\centering{}\includegraphics[width=0.85\textwidth]{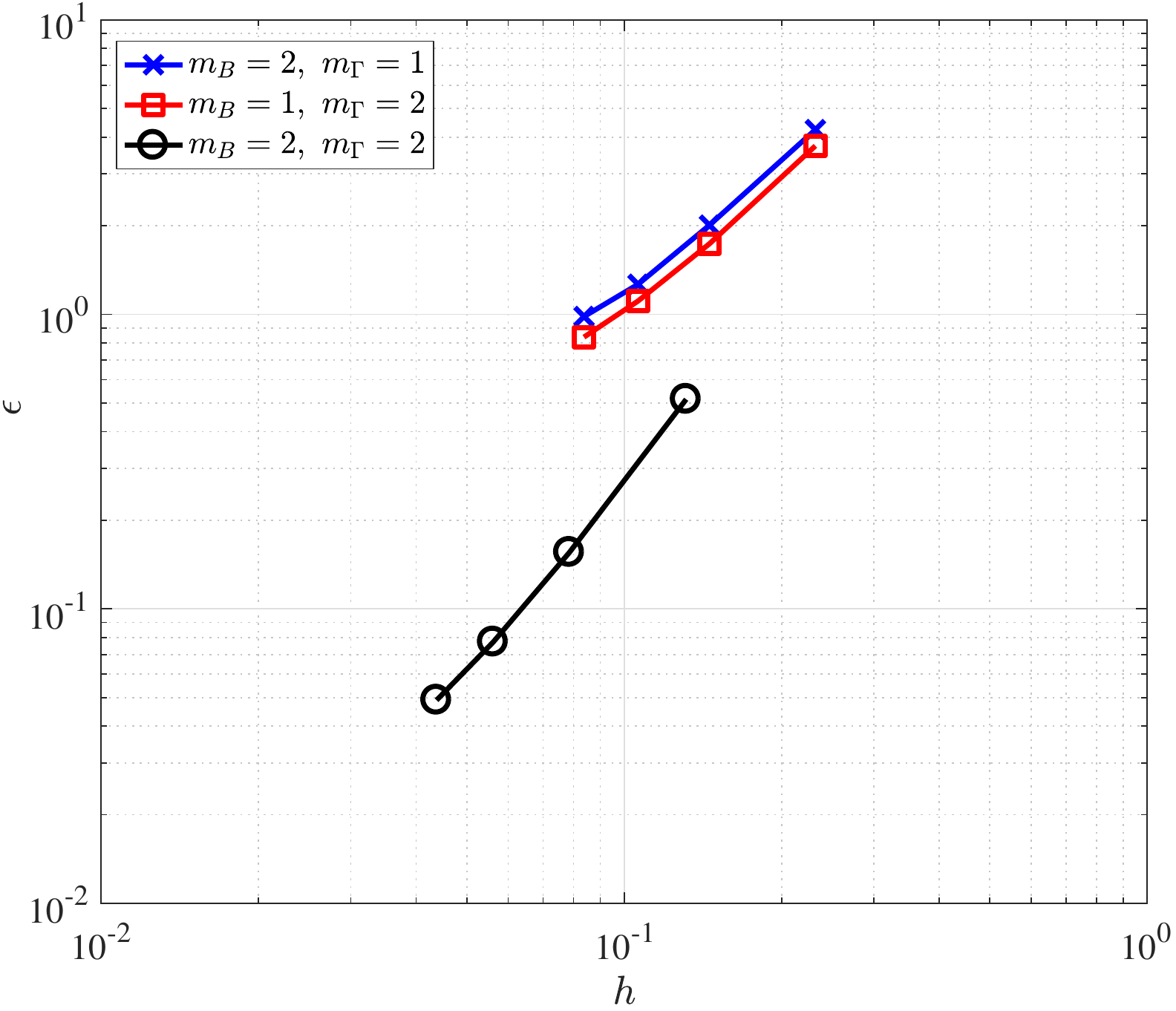}\caption{Stress error convergence for different surface and bulk orders.\label{fig:Stress-error-convergence}}
\end{figure}

\begin{figure}
\begin{centering}
\subfloat[]{\centering{}\includegraphics[width=0.5\textwidth]{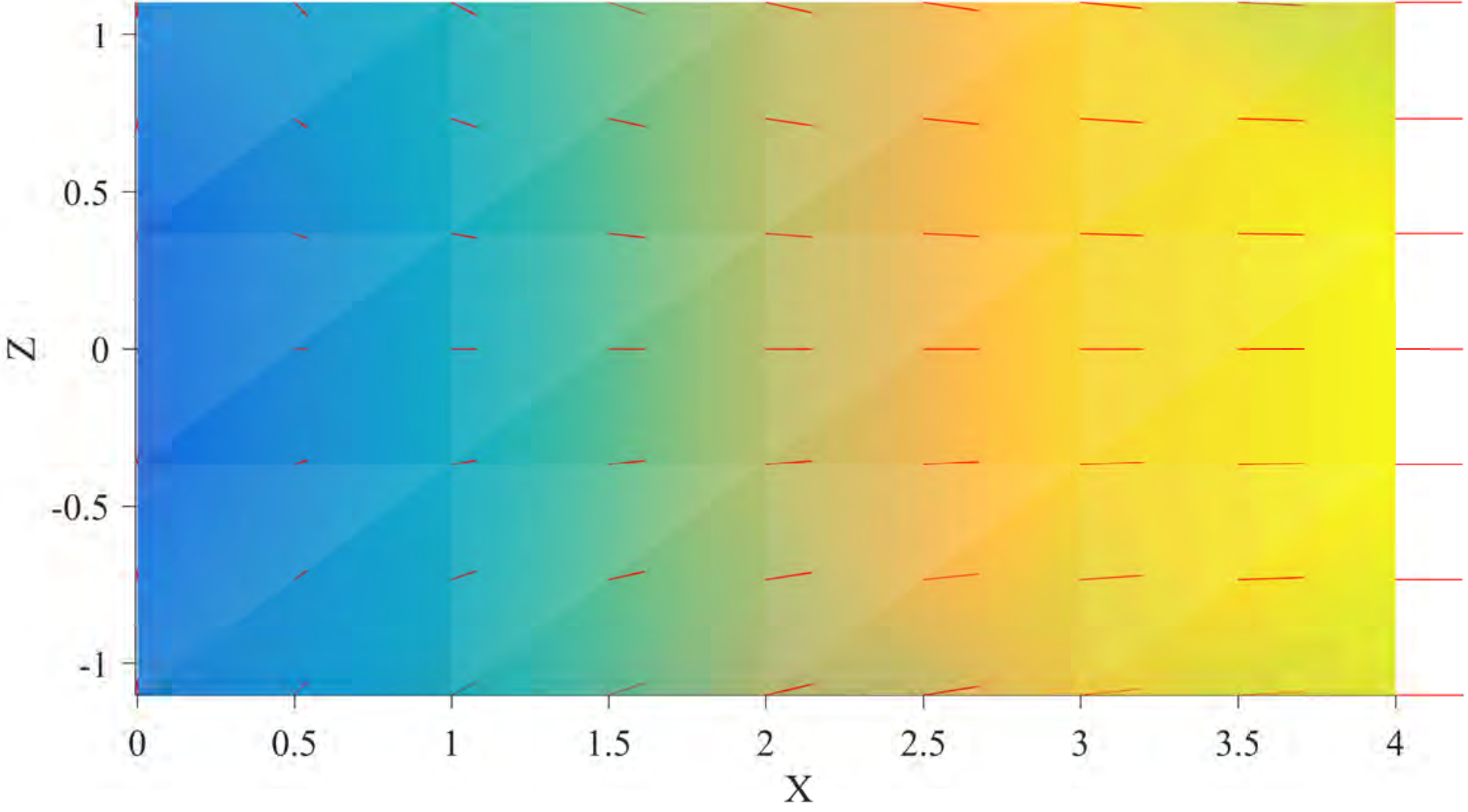}}
\par\end{centering}
\centering{}\subfloat[]{\centering{}\includegraphics[width=0.8\textwidth]{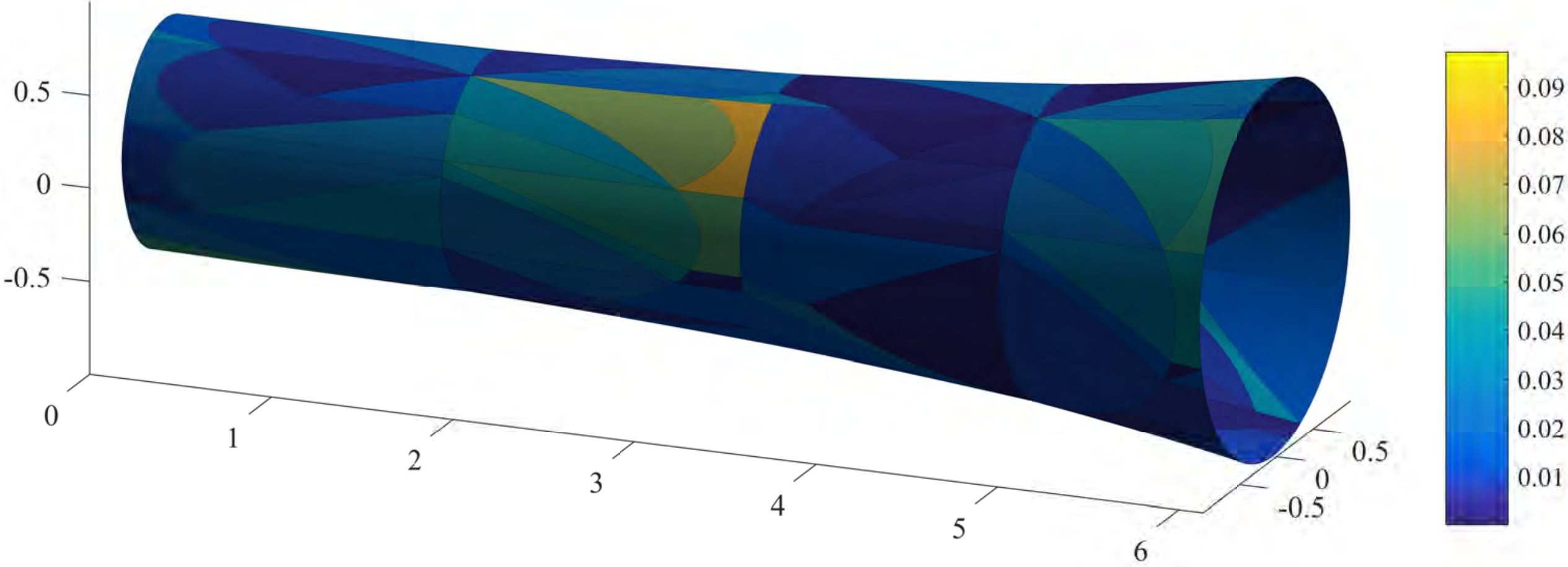}}\caption{Membrane solution with $m_{B}=2$, $m_{\Gamma}=2$. a) On the background
mesh, side view. The color field represents the resultant of the displacements.
b) Solution interpolated to the surface where the color field represents
the stress error $\epsilon_{\sigma}$.\label{fig:Membrane-solution}}
\end{figure}

We use the same example used in \cite{Cenanovic201698,Hansbo2014}.
A cylinder membrane with a radius $r=1$, thickness $t=0.01$ and
length $L=4$, with open ends at $x=0$, $x=L$, with fixed axial
displacements at $x=0$ and radial at $x=L$ and carrying an axial
surface load per unit area

\[
f(x,y,z)=\dfrac{Fx}{2\pi rL^{2}},
\]
where $F=1$ has the unit of force. The material properties are $E=100$,
$\nu=1/2$. The axial stress is given by

\begin{equation}
\sigma_{e}=\dfrac{F\left(1-(x/L)^{2}\right)}{4\pi rt}\label{eq:ExactStress}
\end{equation}
and in the tables and figures $\sigma_{a}$ denotes the approximative
stress computed by

\[
\sigma_{a}:=|\bm{\sigma}_{\Gamma,a}|,\quad\bm{\sigma}_{\Gamma,a}:=[\sigma_{x_{\Gamma}},\sigma_{y_{\Gamma}},\sigma_{z_{\Gamma}}],
\]
where $\sigma_{x_{\Gamma}}$, $\sigma_{y_{\Gamma}}$ and $\sigma_{z_{\Gamma}}$
are the eigenvalues to $\bm{\sigma}_{\Gamma}$. The stress error is
given by

\[
\epsilon_{\sigma}=\|\sigma_{e}-\sigma_{a}\|_{L_{2}(\Gamma_{h})},
\]
see Figure \ref{fig:Stress-error-convergence} for the stress error
convergence. The solution fields using a second order interpolant
can be seen in Figure \ref{fig:Membrane-solution}.

\subsection{Error as a function of the stability factor}

In order to investigate the relation between membrane error and stabilization
factor, we employ two different optimization algorithms. In the first
case where $m_{B}=2$, $m_{\Gamma}=2$ in which we have $j_{h}(\bm{v},\bm{w})=\gamma_{1}j_{h,1}(\bm{v},\bm{w})+\gamma_{2}j_{h,2}(\bm{v},\bm{w})$
the minimization problem is defined as
\[
\mathbb{P}_{1}=\begin{cases}
\min & \epsilon(\gamma_{1},\gamma_{2})\\
\mathrm{s.t.} & 0\leq\gamma_{1}<\infty\\
 & 0\leq\gamma_{2}<\infty
\end{cases}.
\]
The starting point for the minimization using the simplex algorithm
\cite{Nelder1965} is $\bm{\gamma}_{0}=[1,1]$.

\begin{figure}
\centering{}\includegraphics[width=1\textwidth]{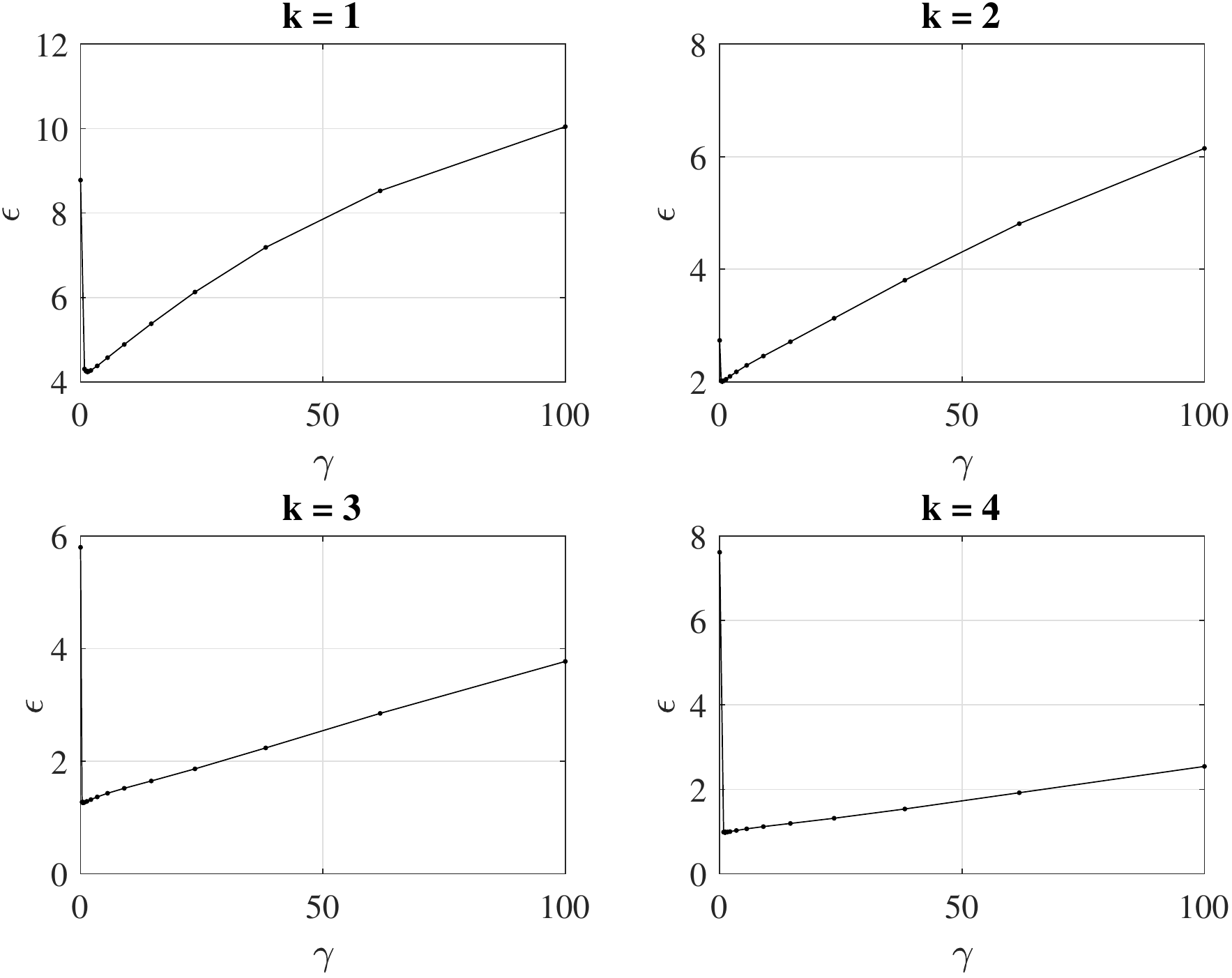}\caption{Stress error as a function of the stabilization factor $\gamma$ for
membrane with $m_{B}=1$, $m_{\Gamma}=1$\label{fig:GammaConvP1P1}}
\end{figure}

\begin{figure}
\centering{}\includegraphics[width=1\textwidth]{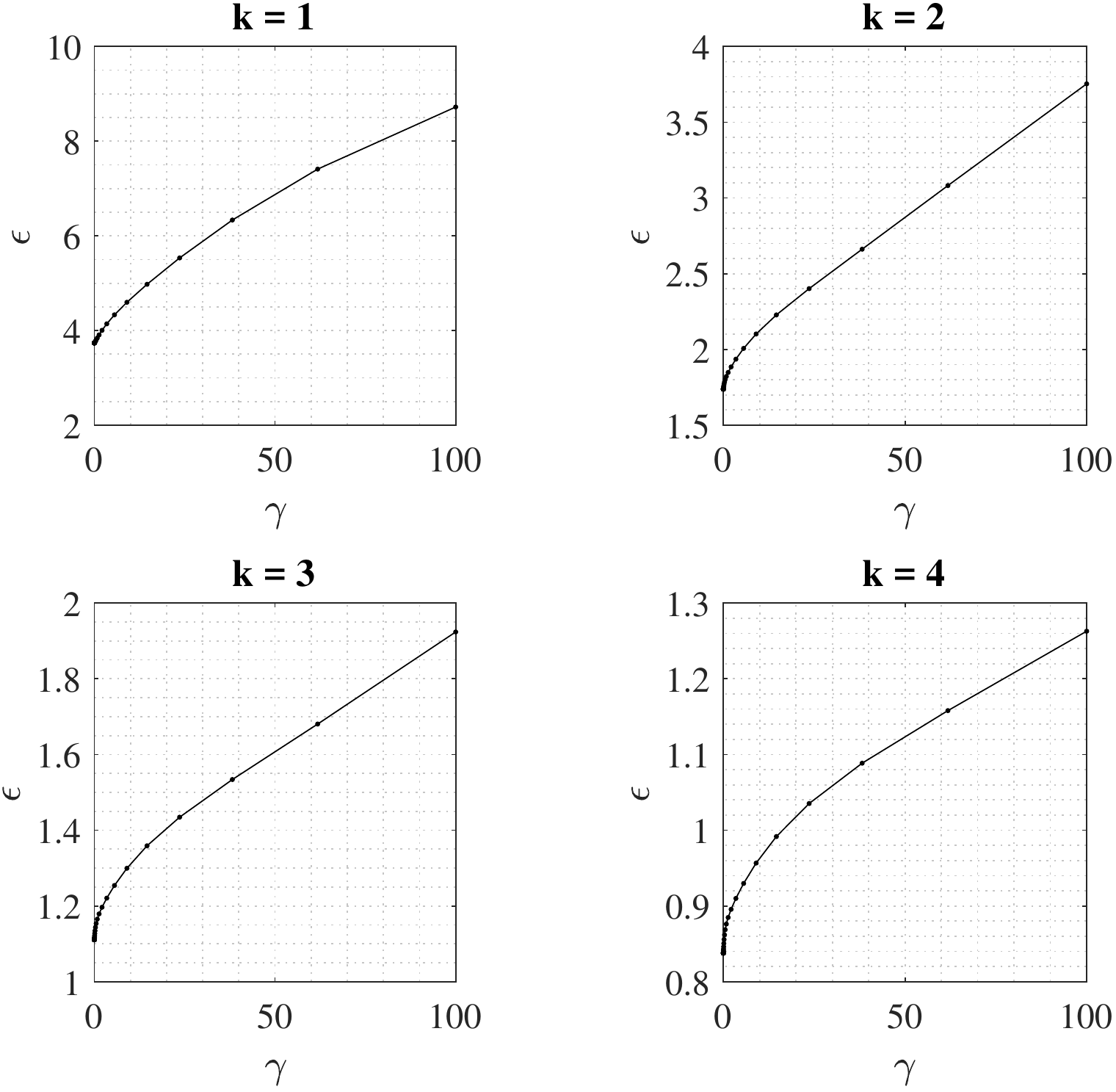}\caption{Stress error as a function of the stabilization factor $\gamma$ for
membrane with $m_{B}=1$, $m_{\Gamma}=2$\label{fig:GammaConvP1P2}}
\end{figure}

\begin{figure}
\begin{centering}
\subfloat[3D view]{\centering{}\includegraphics[width=0.49\textwidth]{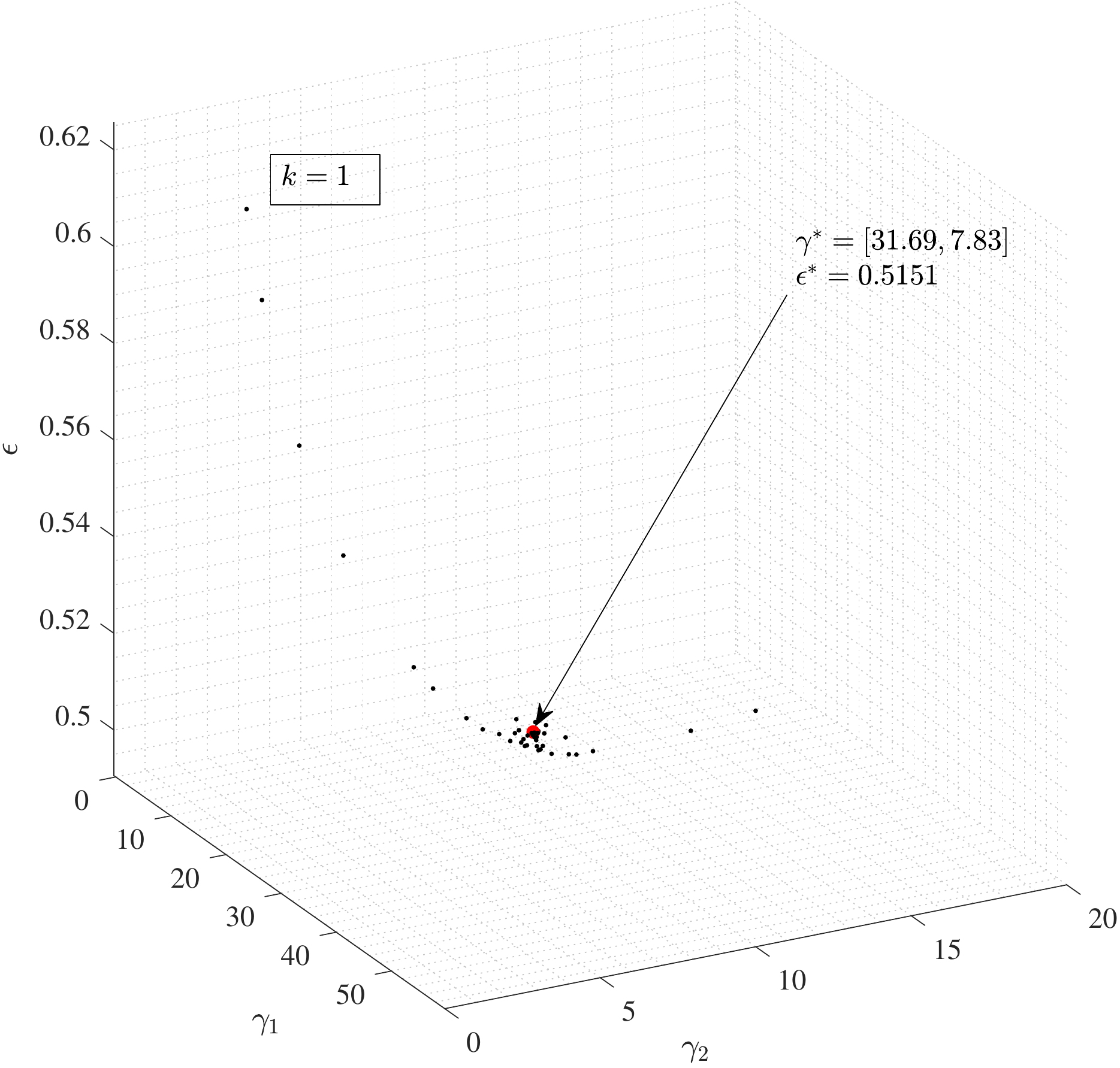}}\subfloat[$\gamma_{2}-\gamma_{1}$ view]{\centering{}\includegraphics[width=0.49\textwidth]{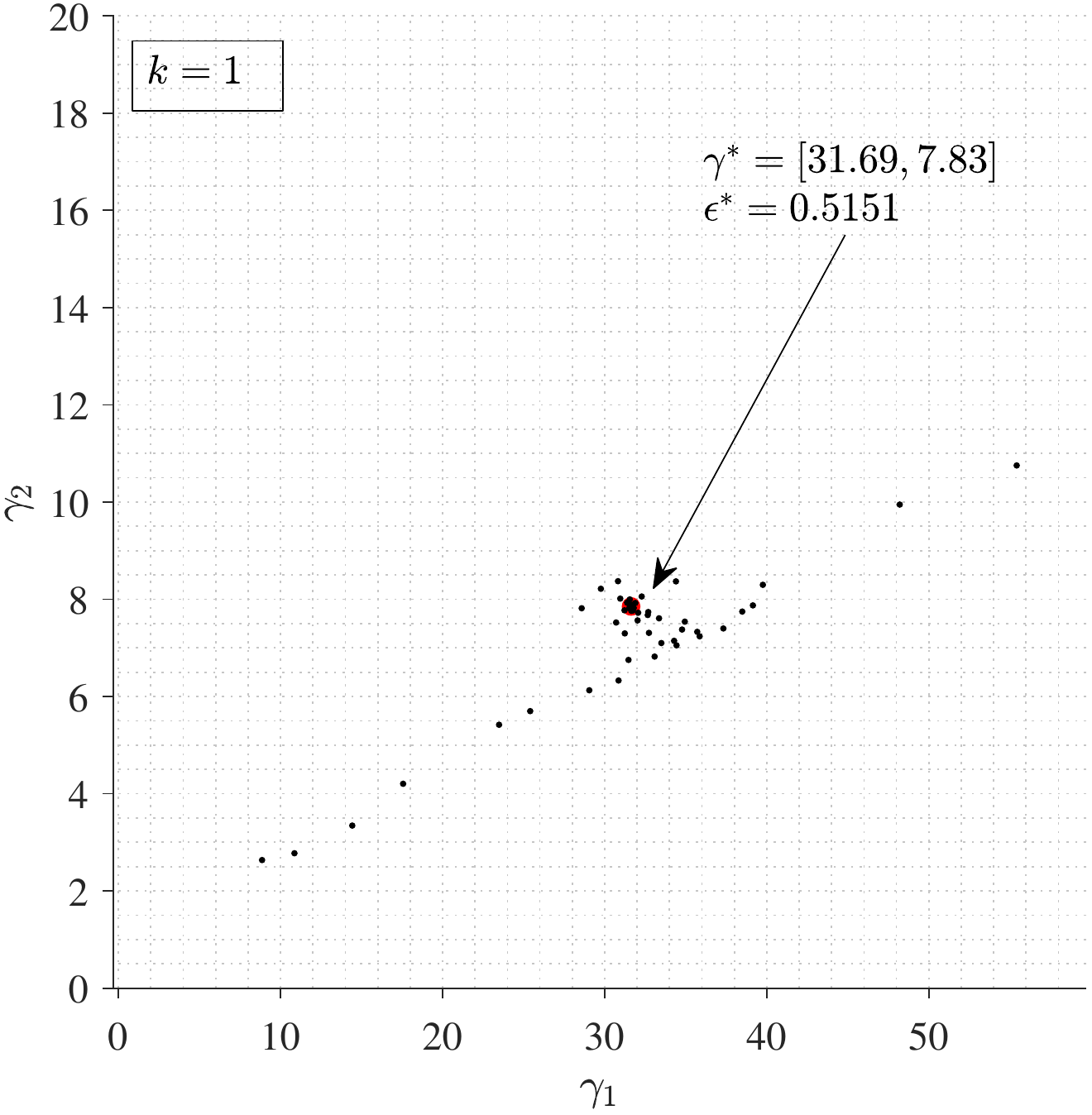}}
\par\end{centering}
\centering{}\subfloat[$\epsilon-\gamma_{1}$ view]{\centering{}\includegraphics[width=0.49\textwidth]{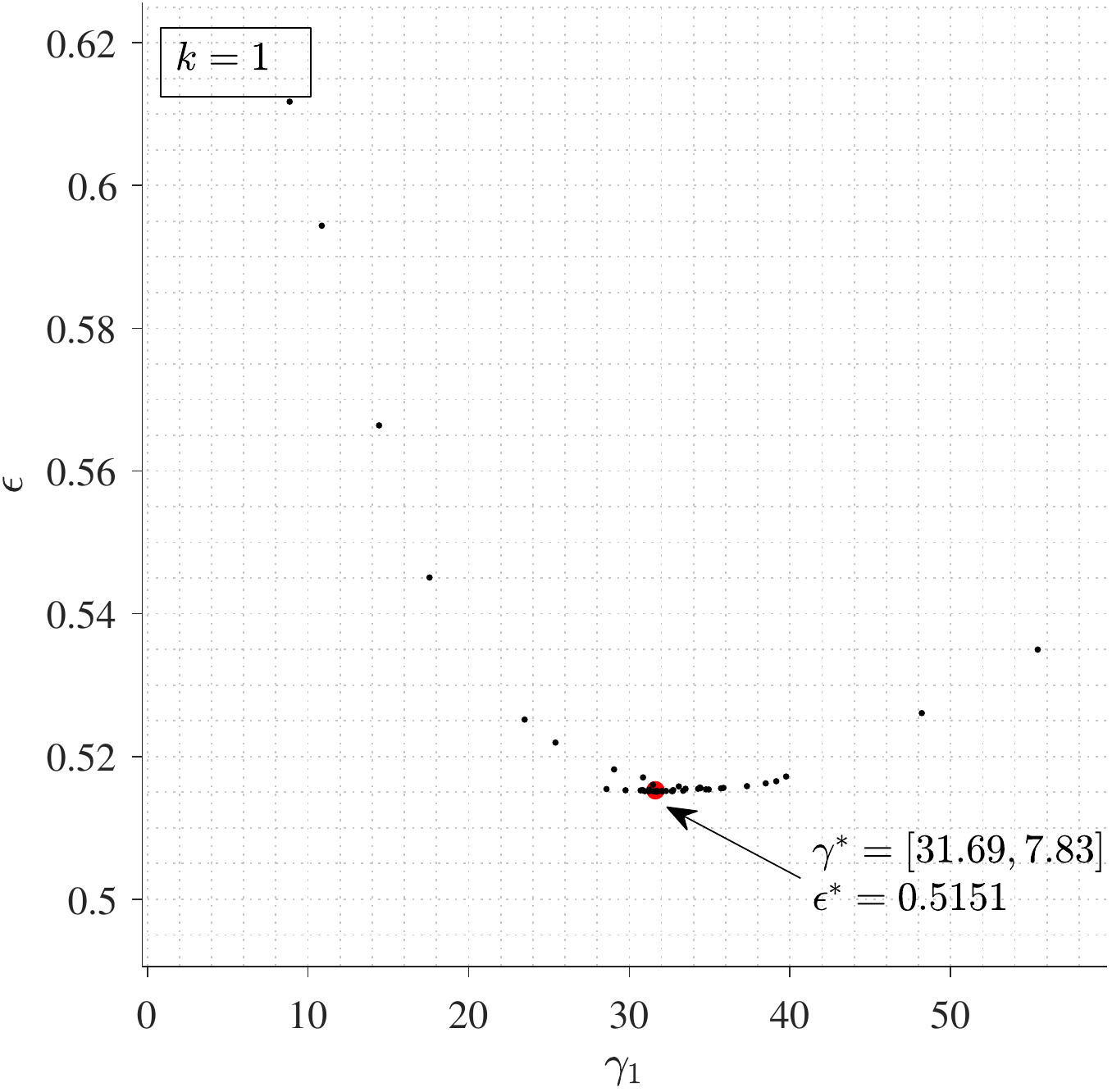}}\subfloat[$\epsilon-\gamma_{2}$ view]{\centering{}\includegraphics[width=0.49\textwidth]{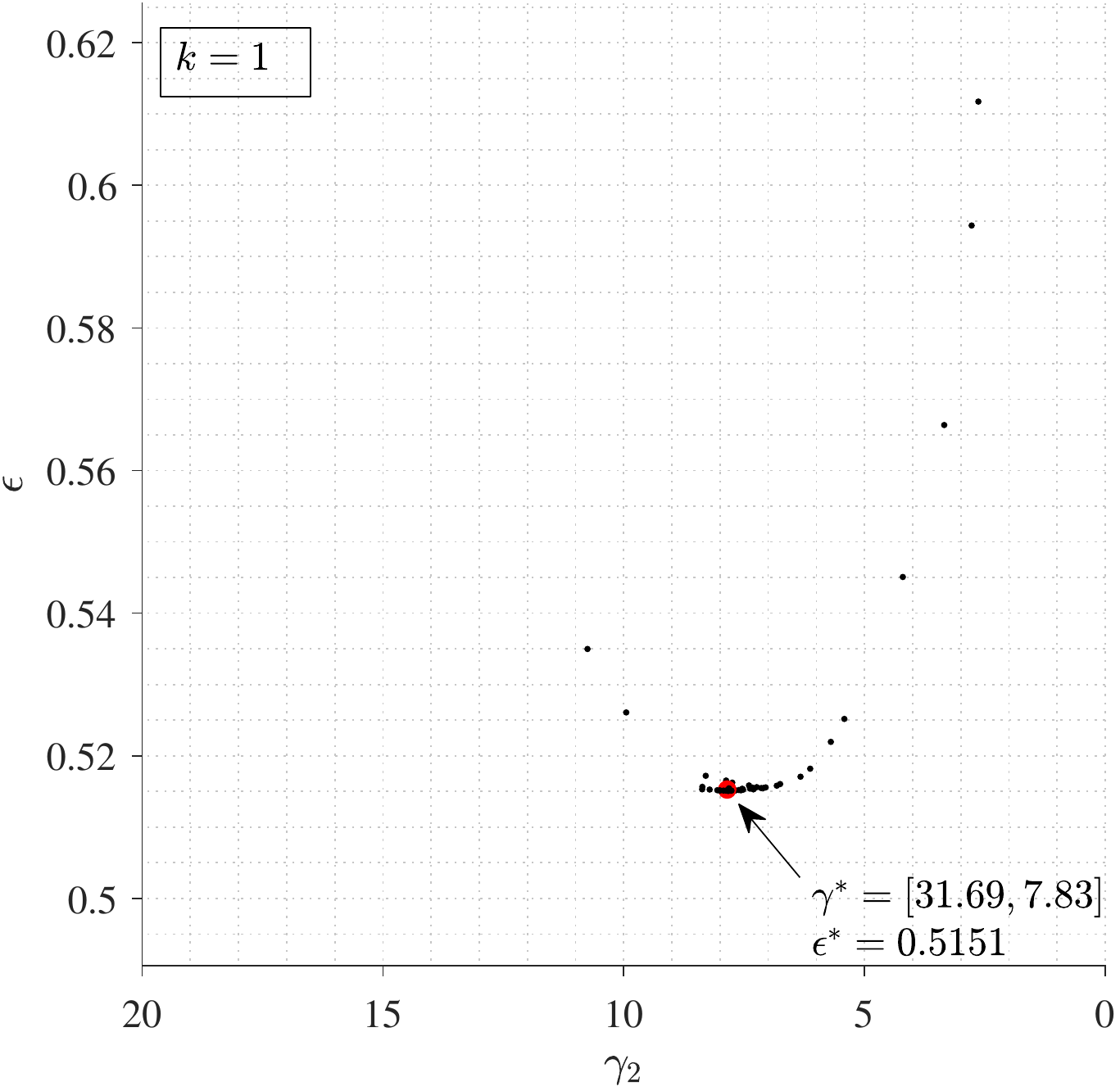}}\caption{Stress error as a function of the stabilization factors $\gamma_{1}$
and $\gamma_{2}$ for membrane with $m_{B}=2$, $m_{\Gamma}=2$ and
$k=1$.\label{fig:GammaConvP1P2K1}}
\end{figure}

\begin{figure}
\begin{centering}
\subfloat[3D view]{\centering{}\includegraphics[width=0.49\textwidth]{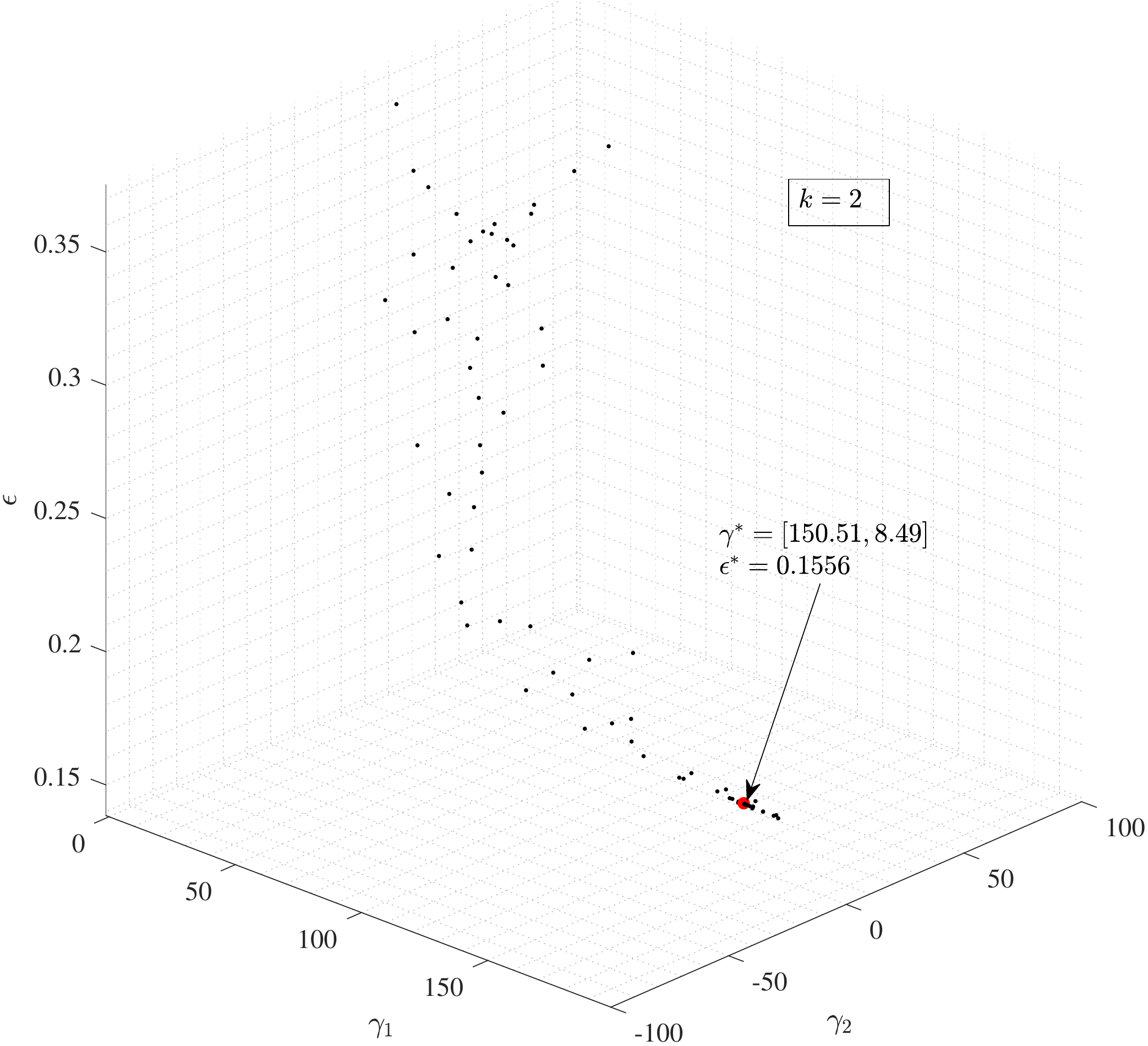}}\subfloat[$\gamma_{2}-\gamma_{1}$ view]{\centering{}\includegraphics[width=0.49\textwidth]{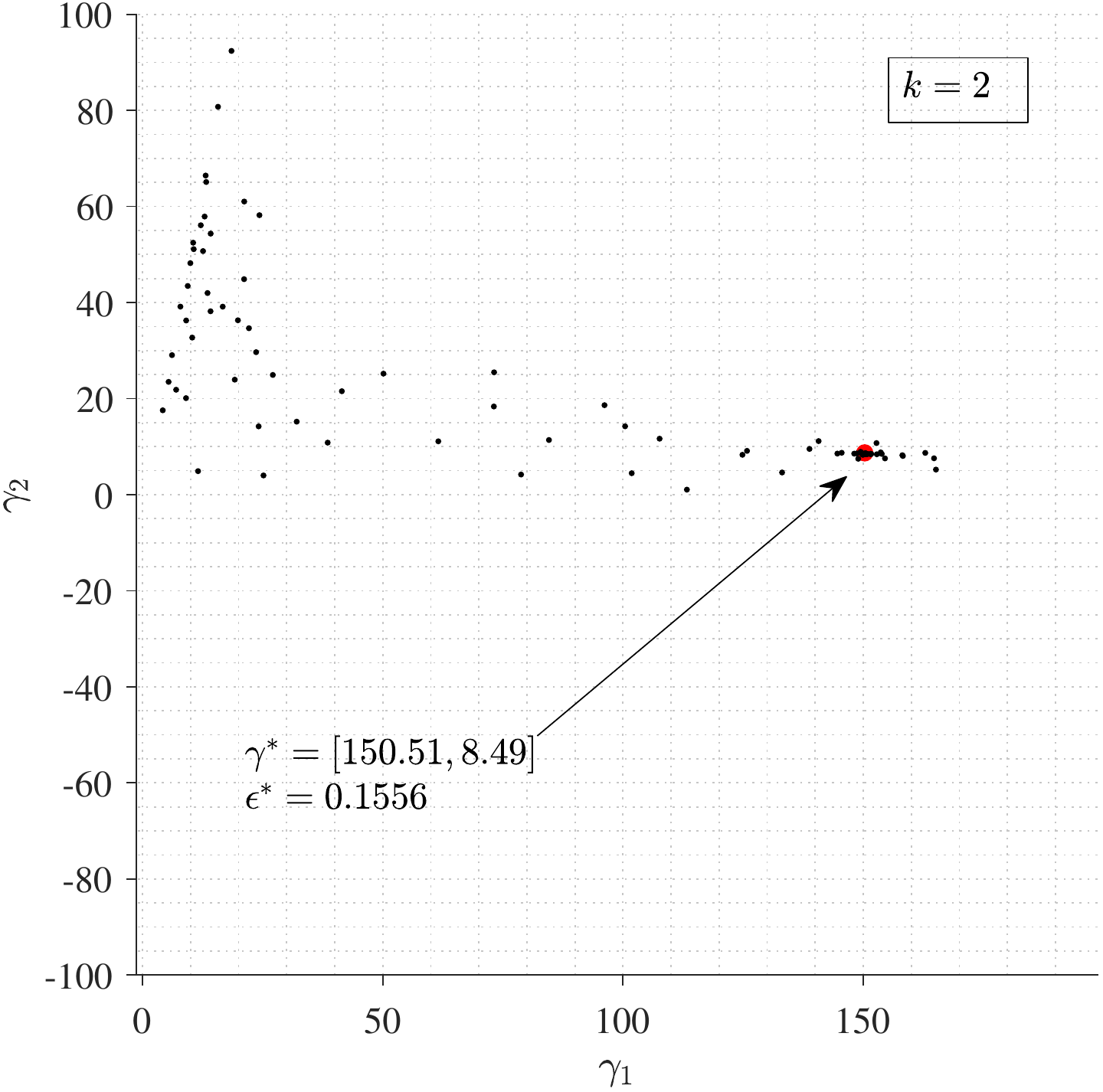}}
\par\end{centering}
\centering{}\subfloat[$\epsilon-\gamma_{1}$ view]{\centering{}\includegraphics[width=0.49\textwidth]{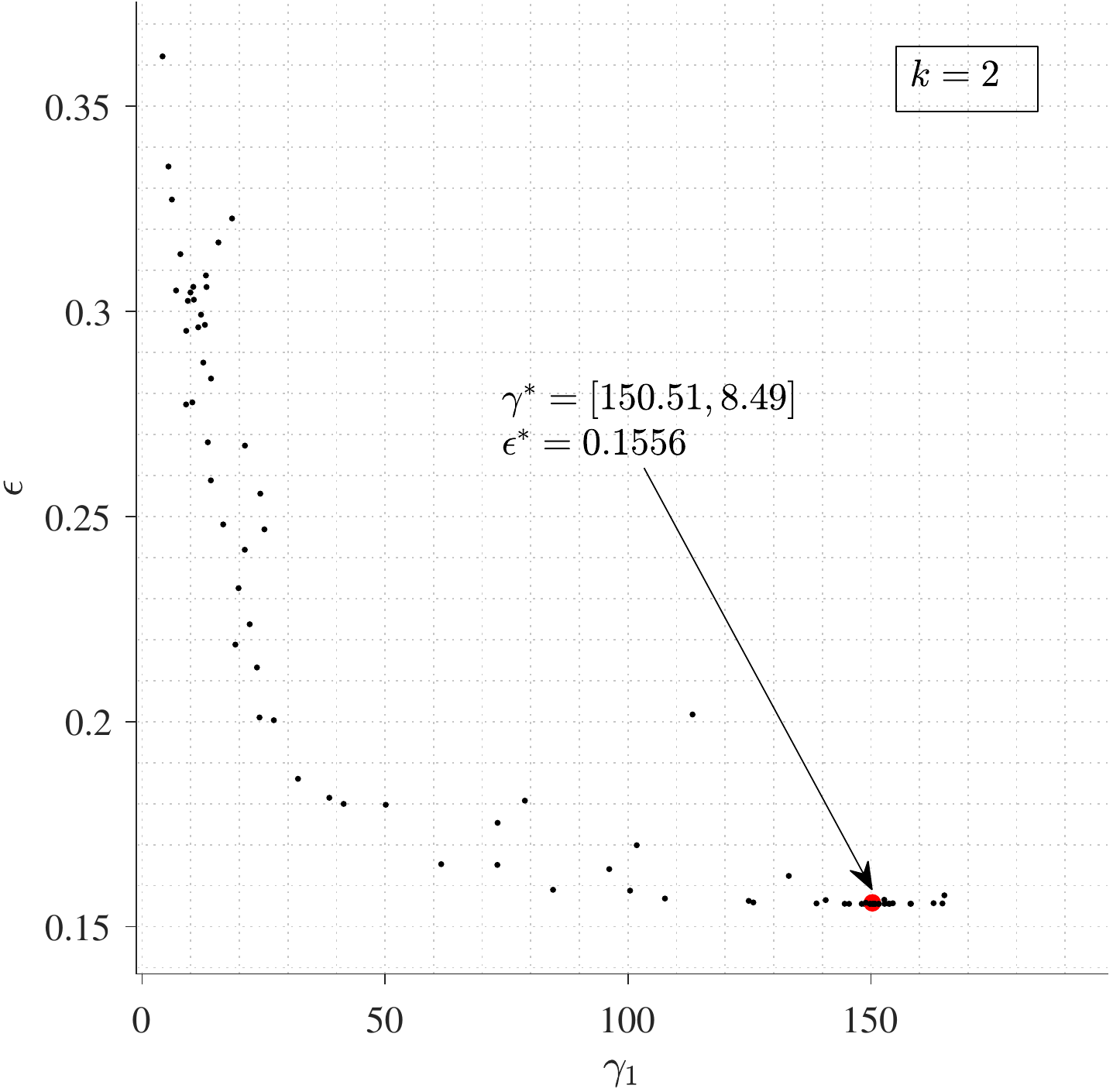}}\subfloat[$\epsilon-\gamma_{2}$ view]{\centering{}\includegraphics[width=0.49\textwidth]{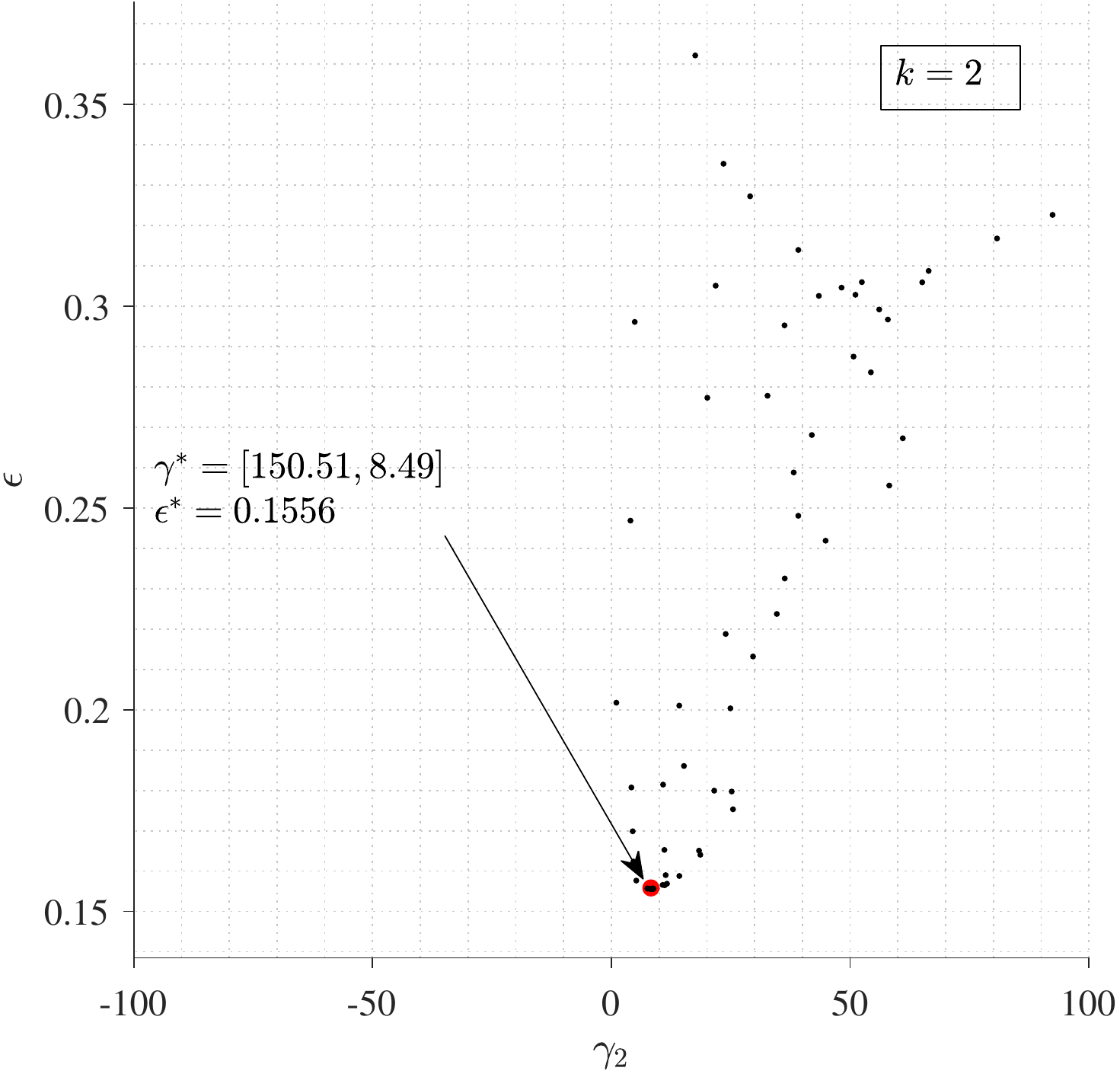}}\caption{Stress error as a function of the stabilization factors $\gamma_{1}$
and $\gamma_{2}$ for membrane with $m_{B}=2$, $m_{\Gamma}=2$ and
$k=2$.\label{fig:GammaConvP1P2K2}}
\end{figure}

\begin{figure}
\begin{centering}
\subfloat[3D view]{\centering{}\includegraphics[width=0.49\textwidth]{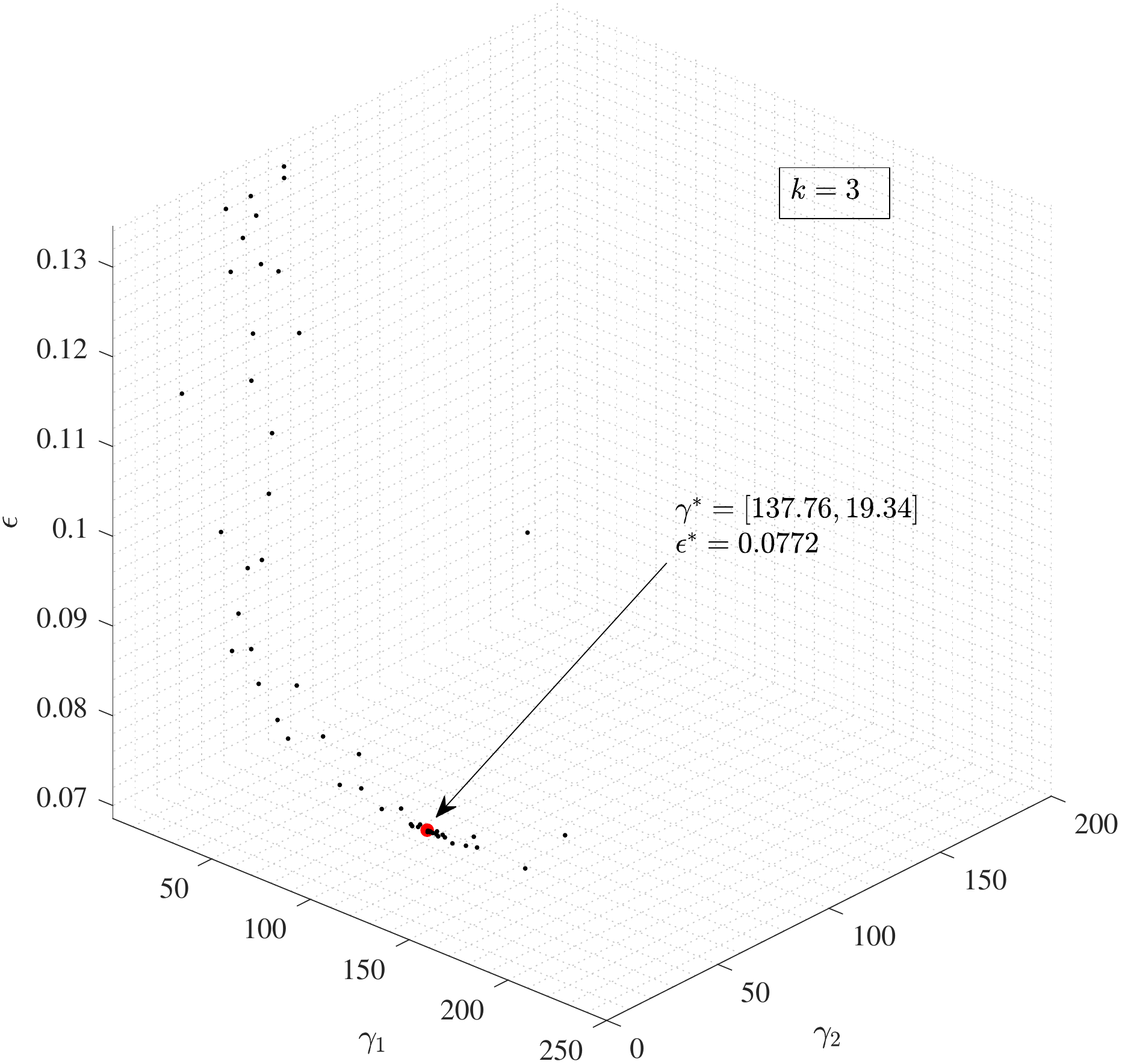}}\subfloat[$\gamma_{2}-\gamma_{1}$ view]{\centering{}\includegraphics[width=0.49\textwidth]{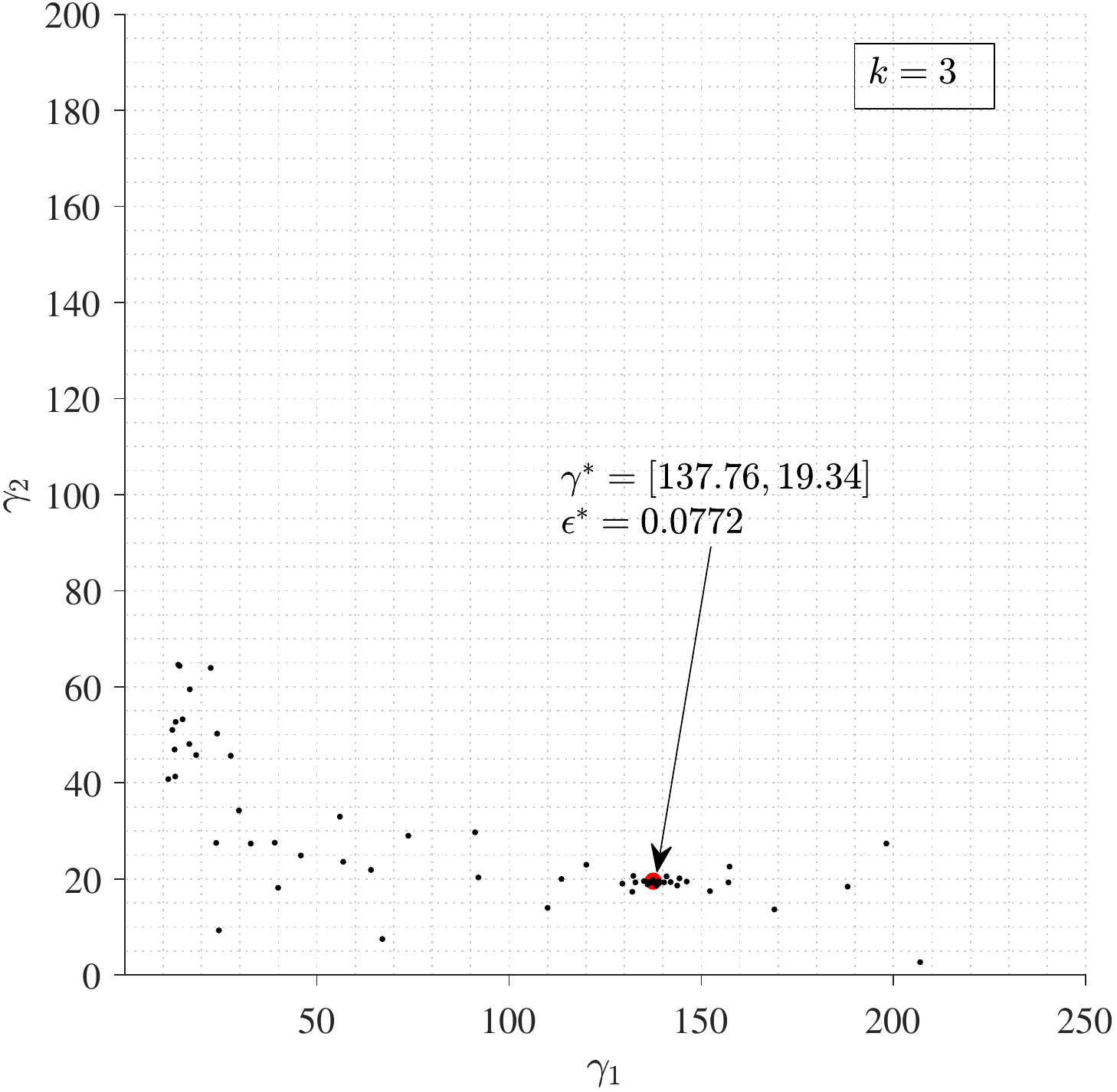}}
\par\end{centering}
\centering{}\subfloat[$\epsilon-\gamma_{1}$ view]{\centering{}\includegraphics[width=0.49\textwidth]{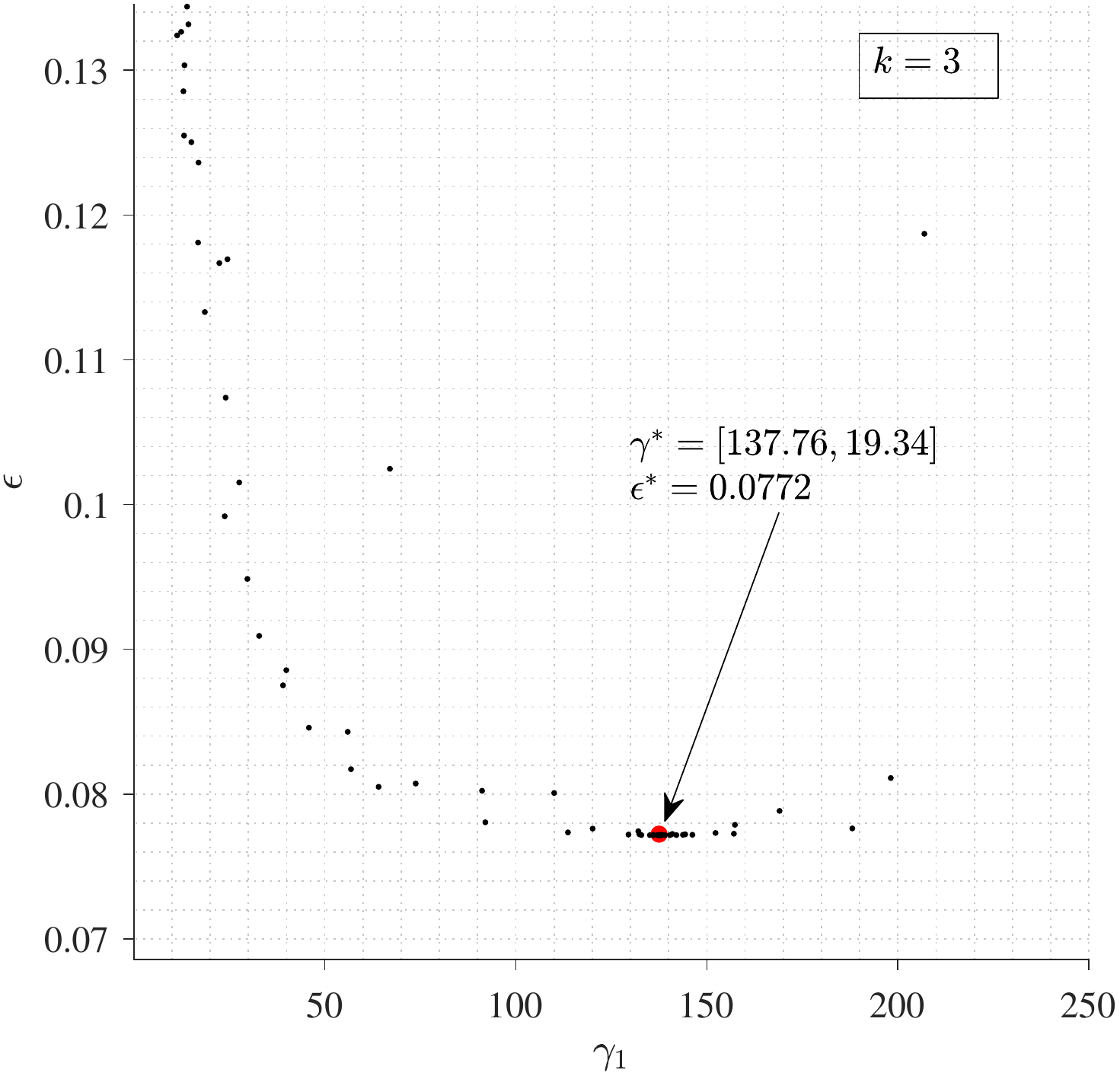}}\subfloat[$\epsilon-\gamma_{2}$ view]{\centering{}\includegraphics[width=0.49\textwidth]{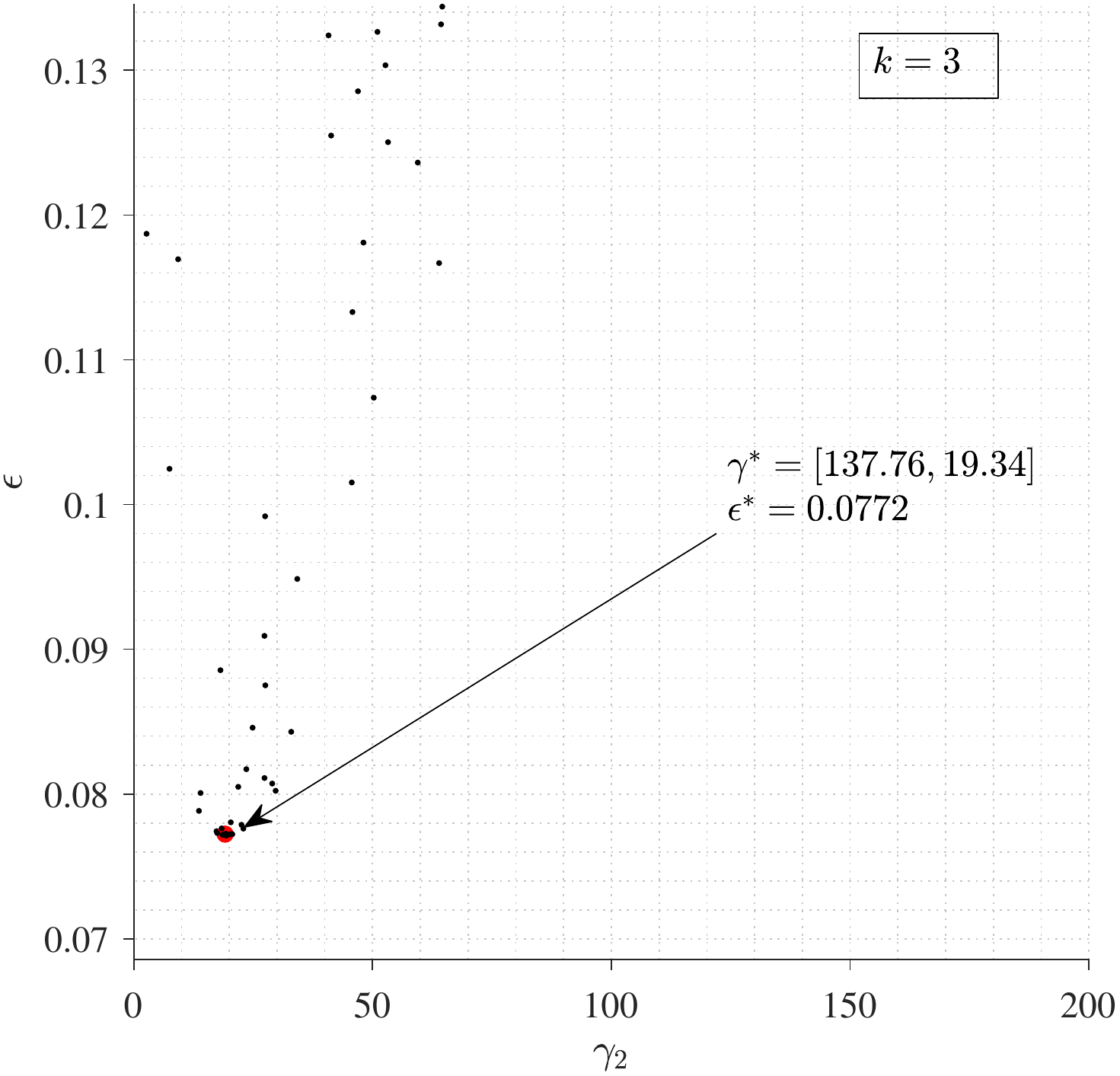}}\caption{Stress error as a function of the stabilization factors $\gamma_{1}$
and $\gamma_{2}$ for membrane with $m_{B}=2$, $m_{\Gamma}=2$ and
$k=3$.\label{fig:GammaConvP1P2K3}}
\end{figure}

\begin{figure}
\begin{centering}
\subfloat[3D view]{\centering{}\includegraphics[width=0.49\textwidth]{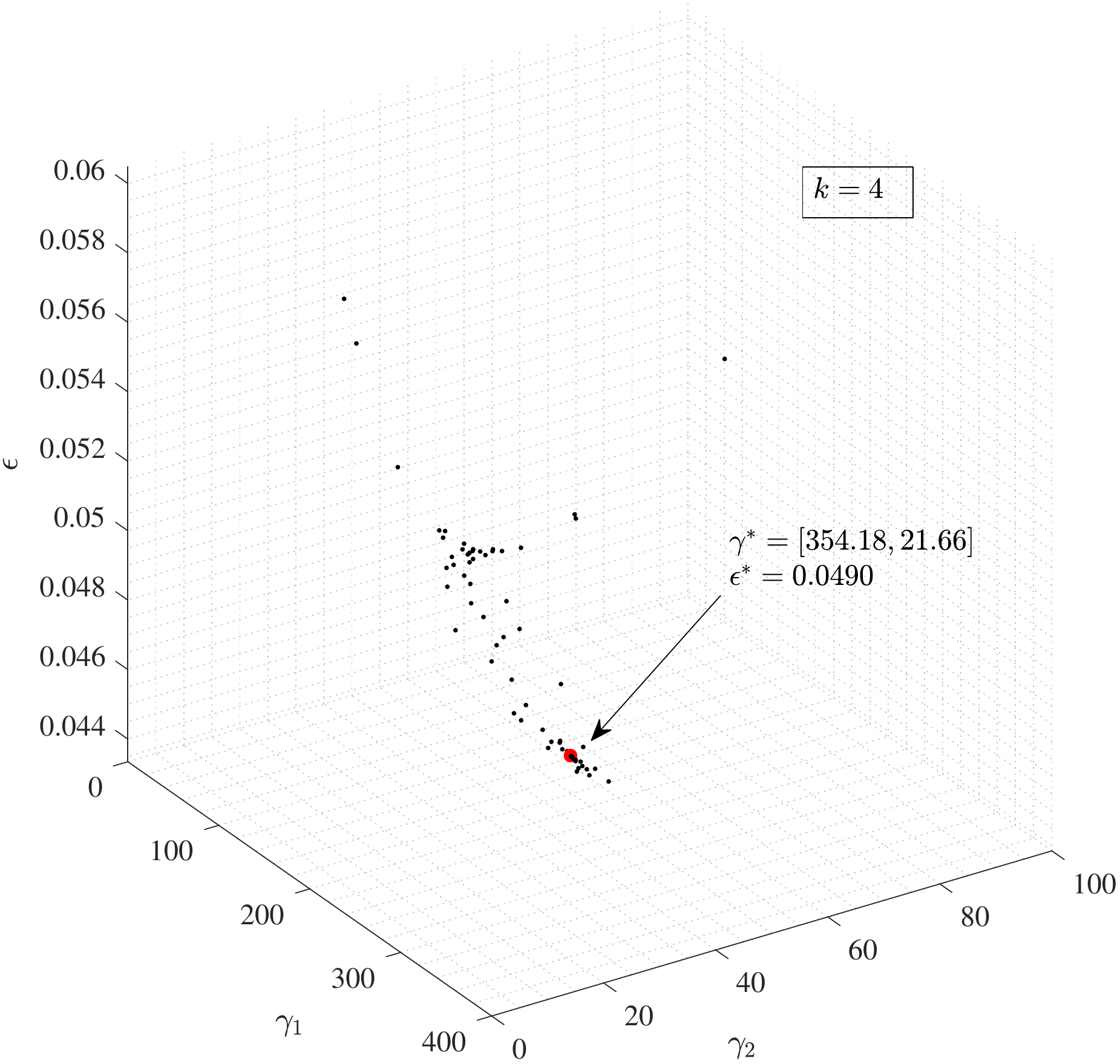}}\subfloat[$\gamma_{2}-\gamma_{1}$ view]{\centering{}\includegraphics[width=0.49\textwidth]{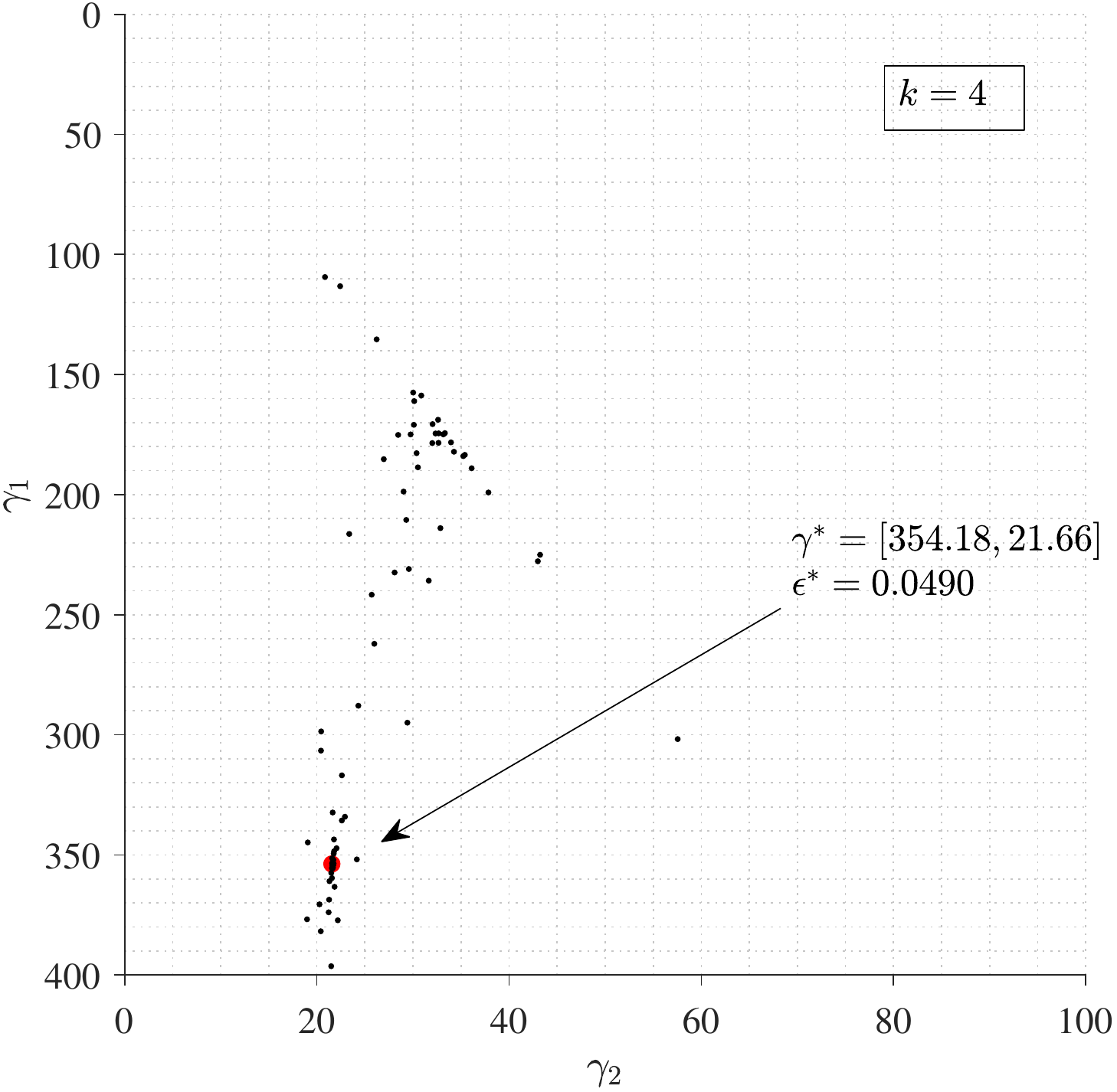}}
\par\end{centering}
\centering{}\subfloat[$\epsilon-\gamma_{1}$ view]{\centering{}\includegraphics[width=0.49\textwidth]{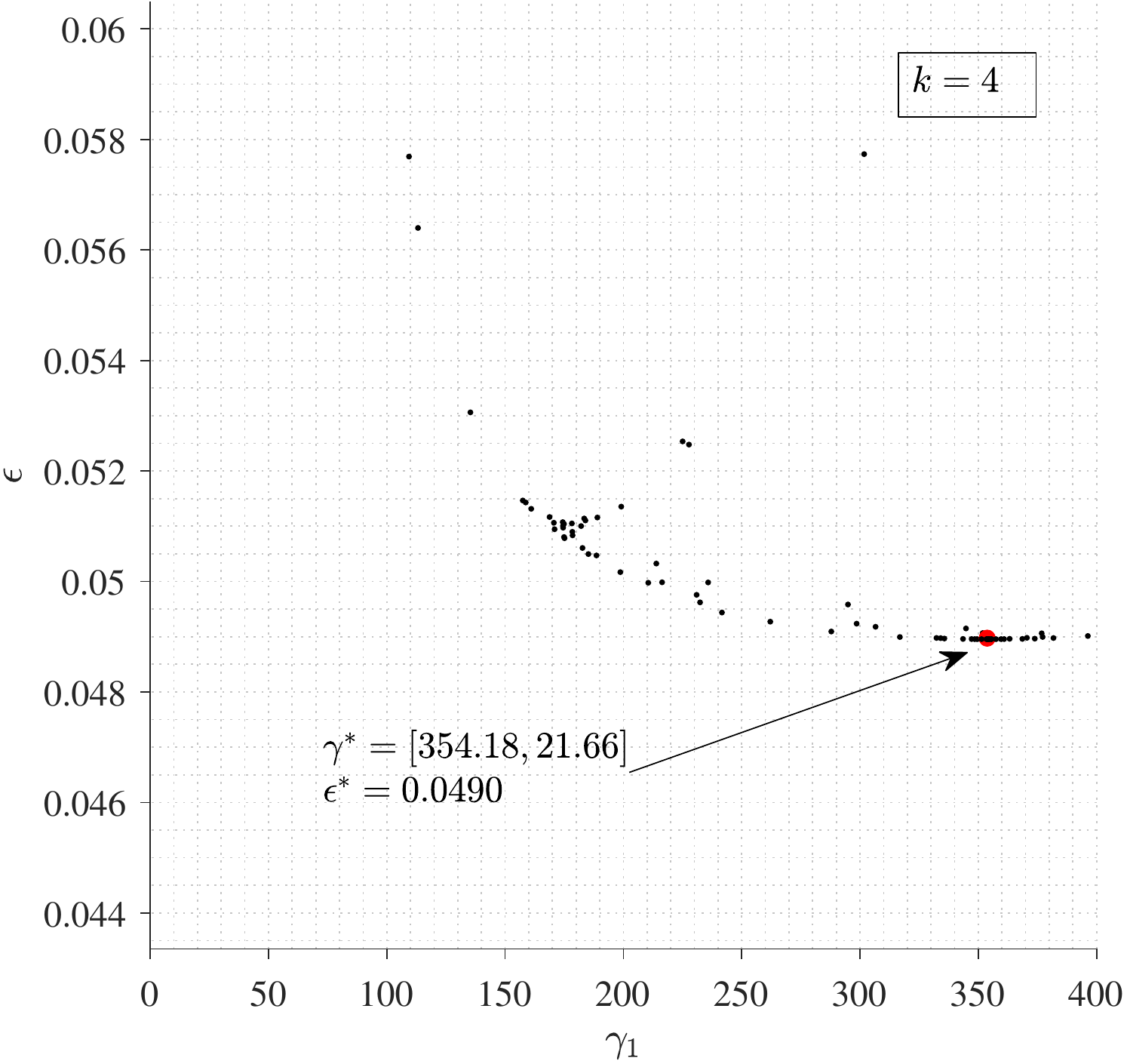}}\subfloat[$\epsilon-\gamma_{2}$ view]{\centering{}\includegraphics[width=0.49\textwidth]{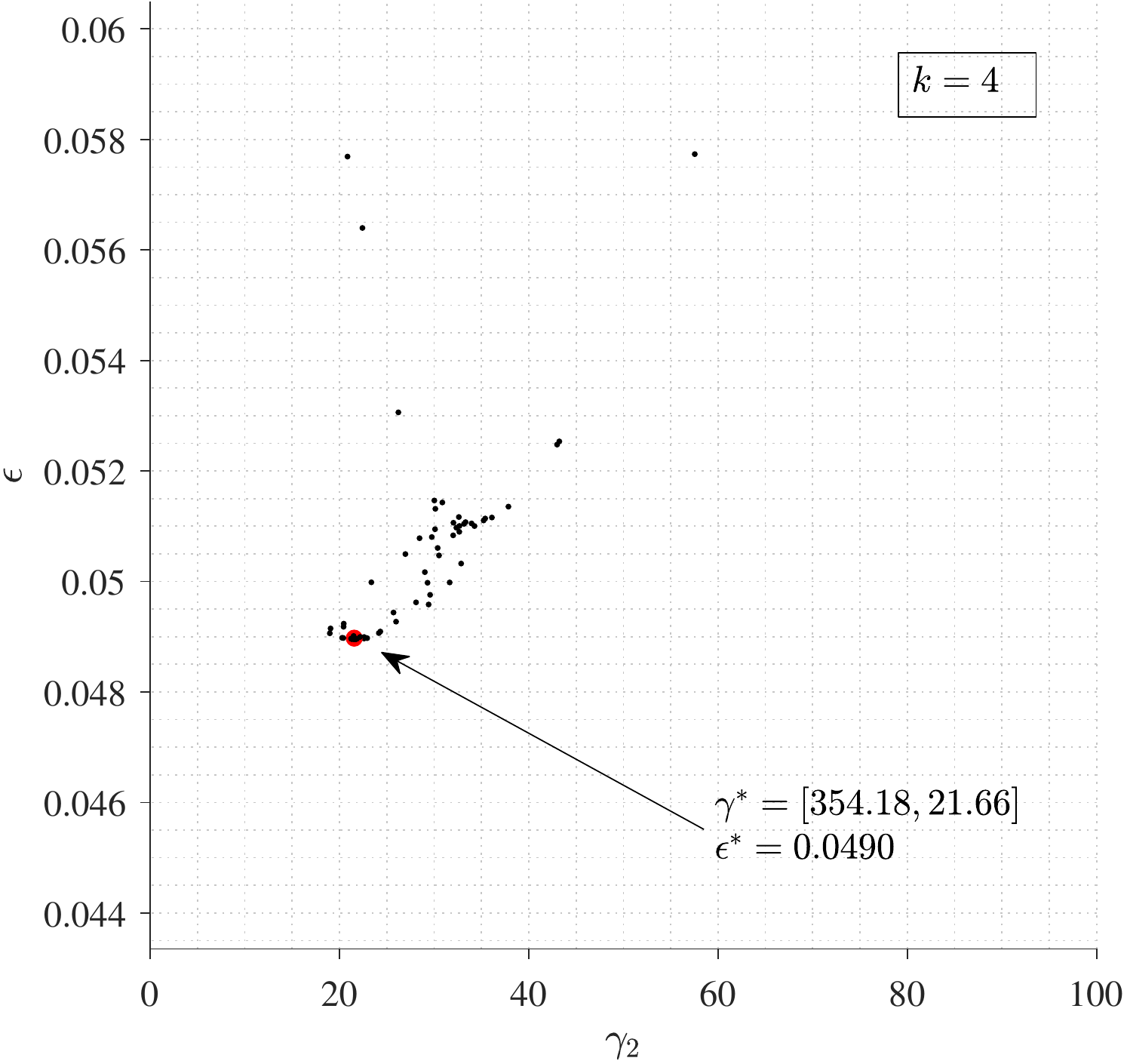}}\caption{Stress error as a function of the stabilization factors $\gamma_{1}$
and $\gamma_{2}$ for membrane with $m_{B}=2$, $m_{\Gamma}=2$ and
$k=4$.\label{fig:GammaConvP1P2K4}}
\end{figure}

\begin{figure}
\centering{}\subfloat[]{\centering{}\includegraphics[width=0.49\textwidth]{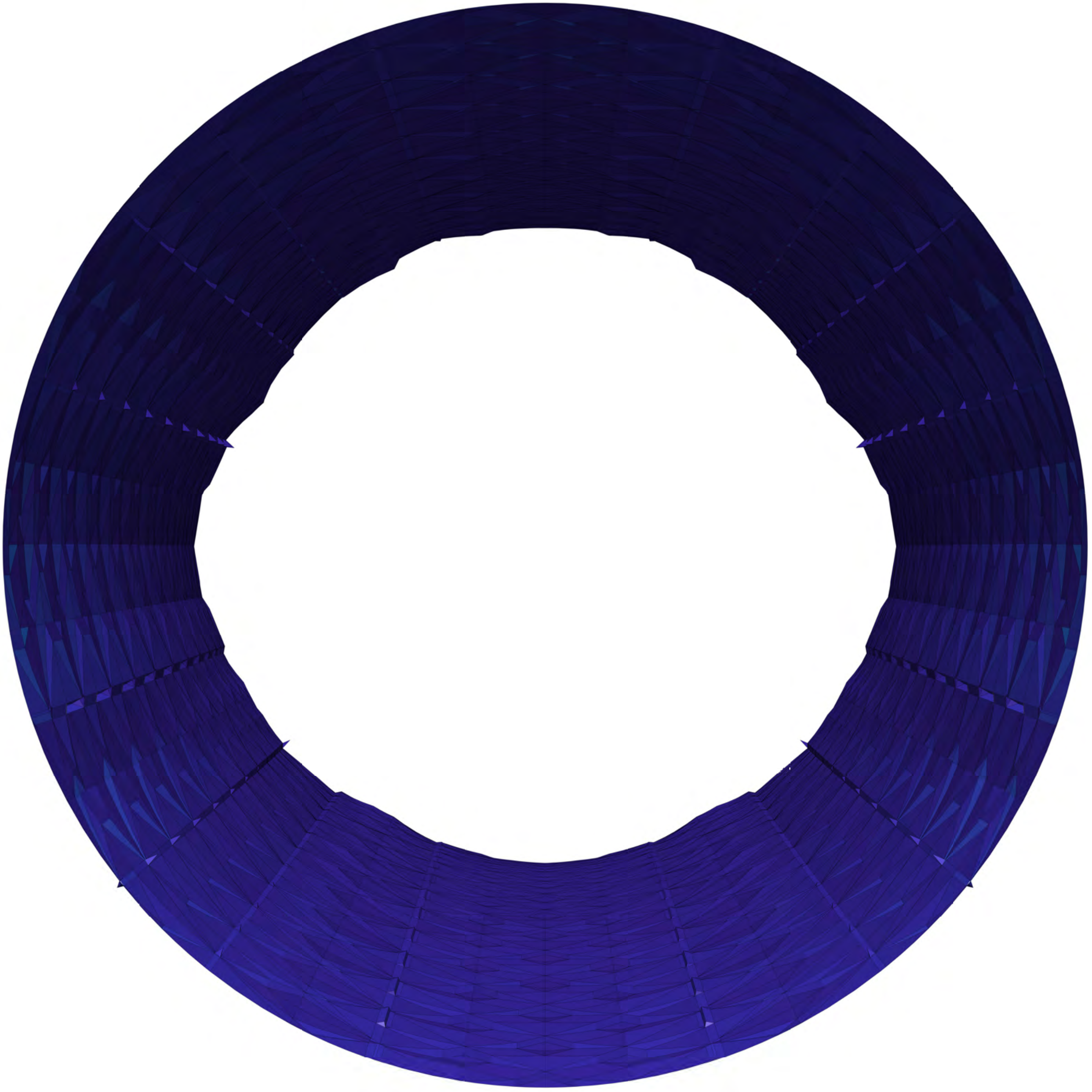}}\hspace{0.01\textwidth}\subfloat[]{\centering{}\includegraphics[width=0.49\textwidth]{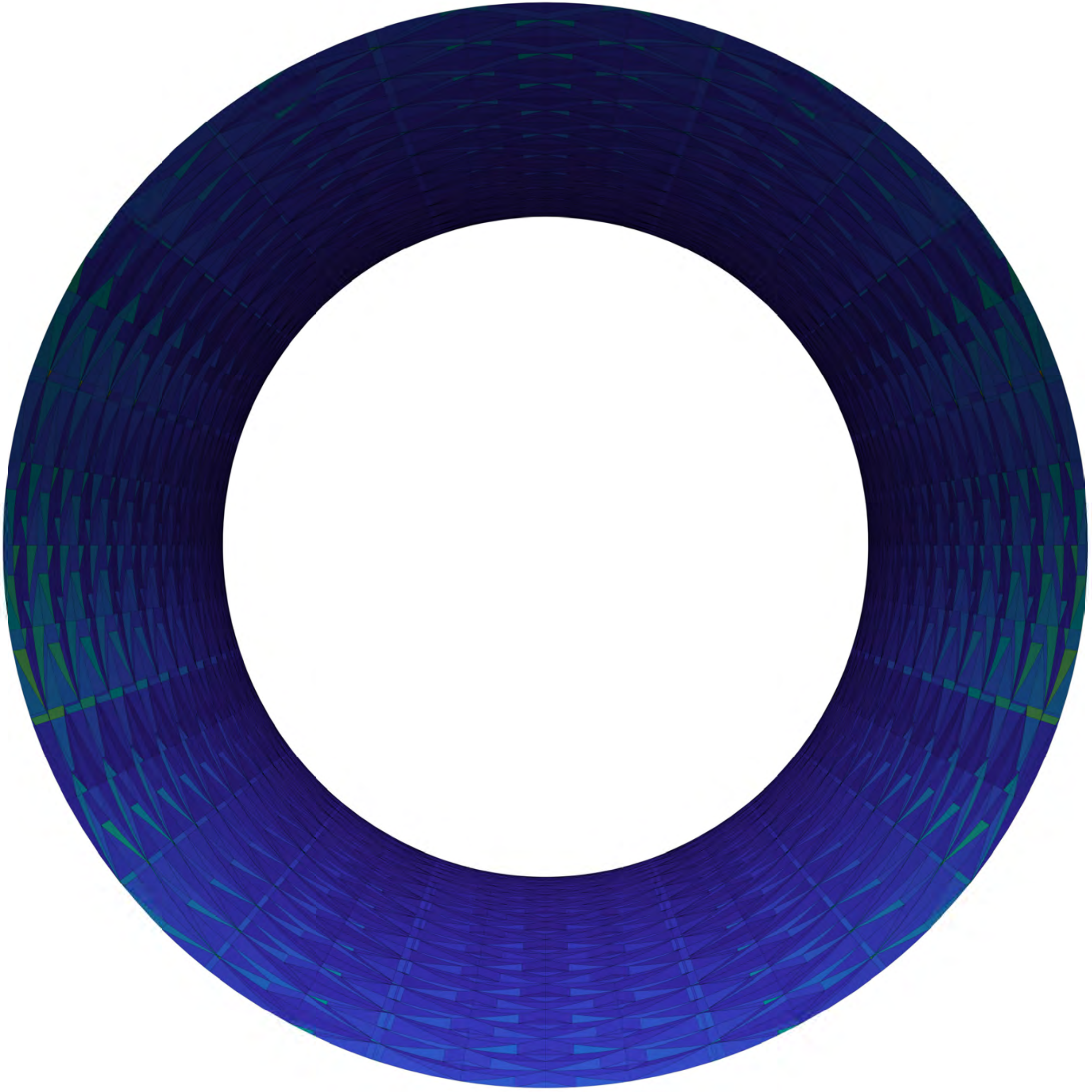}}\caption{Displacement field (10 times enlarged) for membrane with $m_{B}=1$,
$m_{\Gamma}=2$ with. (a) Front view of interpolated displacements
with $\gamma=0$. (b) Front view of interpolated displacements with
$\gamma=10$. \label{fig:Displacement-fields}}
\end{figure}

In case of $m_{B}=1$, $m_{\Gamma}=1$ and $m_{B}=1$, $m_{\Gamma}=2$
we have $j_{h}(\bm{v},\bm{w})=\gamma j_{h}(\bm{v},\bm{w})$ and define
the minimization problem as

\[
\mathbb{P}_{2}=\begin{cases}
\min & \epsilon(\gamma)\\
\mathrm{s.t.} & 0\leq\gamma\leq100
\end{cases}.
\]
This problem is solved using the golden search method. In both optimization
problems, the optimal parameter is denoted with the superscript $*$.
The results of this study can be seen in Table \ref{tab:StressErrorConvP1P1}
to Table \ref{tab:StressErrorConvP2P2} and Figure \ref{fig:GammaConvP1P1}
to Figure \ref{fig:GammaConvP1P2K4}. Note that although the solution
with $m_{B}=1$, $m_{\Gamma}=2$ is stable without any stabilization,
the interpolation of $U$ onto $\Gamma_{h}$ is not, see Figure \ref{fig:Displacement-fields}.

\begin{table}
\begin{centering}
\begin{tabular}{|c|c|c|c|c|}
\hline 
$k$ & $h$ & $\epsilon_{\sigma}$ & Rate & $\gamma_{1}^{*}$\tabularnewline
\hline 
\hline 
1 & 0.2321 & 4.2421 & - & 1.4332\tabularnewline
\hline 
2 & 0.1456 & 2.0101 & 1.6017 & 0.5107\tabularnewline
\hline 
3 & 0.1063 & 1.2655 & 1.4708 & 0.5440\tabularnewline
\hline 
4 & 0.0838 & 0.9838 & 1.0587 & 1.0801\tabularnewline
\hline 
\end{tabular}
\par\end{centering}
\centering{}\caption{Error convergence for the membrane with $m_{B}=1$, $m_{\Gamma}=1$\label{tab:StressErrorConvP1P1}}
\end{table}

\begin{table}
\begin{centering}
\begin{tabular}{|c|c|c|c|c|}
\hline 
$k$ & $h$ & $\epsilon_{\sigma}$ & Rate & $\gamma_{1}^{*}$\tabularnewline
\hline 
\hline 
1 & 0.2321 & 3.7366 & - & 0\tabularnewline
\hline 
2 & 0.1456 & 1.7383 & 1.6411 & 0\tabularnewline
\hline 
3 & 0.1063 & 1.1108 & 1.4235 & 0\tabularnewline
\hline 
4 & 0.0838 & 0.8377  & 1.1864 & 0\tabularnewline
\hline 
\end{tabular}
\par\end{centering}
\centering{}\caption{Error convergence for the membrane with $m_{B}=1$, $m_{\Gamma}=2$\label{tab:StressErrorConvP1P2}}
\end{table}

\begin{table}
\begin{centering}
\begin{tabular}{|c|c|c|c|c|c|}
\hline 
$k$ & $h$ & $\epsilon_{\sigma}$ & $\gamma_{1}^{*}$ & $\gamma_{2}^{*}$ & Rate\tabularnewline
\hline 
\hline 
1 & 0.1314 & 0.5151 & 31.6944 & 7.8296 & -\tabularnewline
\hline 
2 & 0.0786 & 0.1556 & 150.5121 & 8.4932 & 2.3295\tabularnewline
\hline 
3 & 0.0562 & 0.0772 & 137.7599 & 19.3374 & 2.0894\tabularnewline
\hline 
4 & 0.0438 & 0.0490 & 354.1755 & 21.6636 & 1.8235\tabularnewline
\hline 
\end{tabular}
\par\end{centering}
\centering{}\caption{$\epsilon(\gamma_{1},\gamma_{2})$ for membrane with $m_{B}=2$, $m_{\Gamma}=2$\label{tab:StressErrorConvP2P2}}
\end{table}

\subsection{Geometrical error}

\begin{figure}
\centering{}\subfloat[]{\begin{centering}
\includegraphics[width=0.49\textwidth]{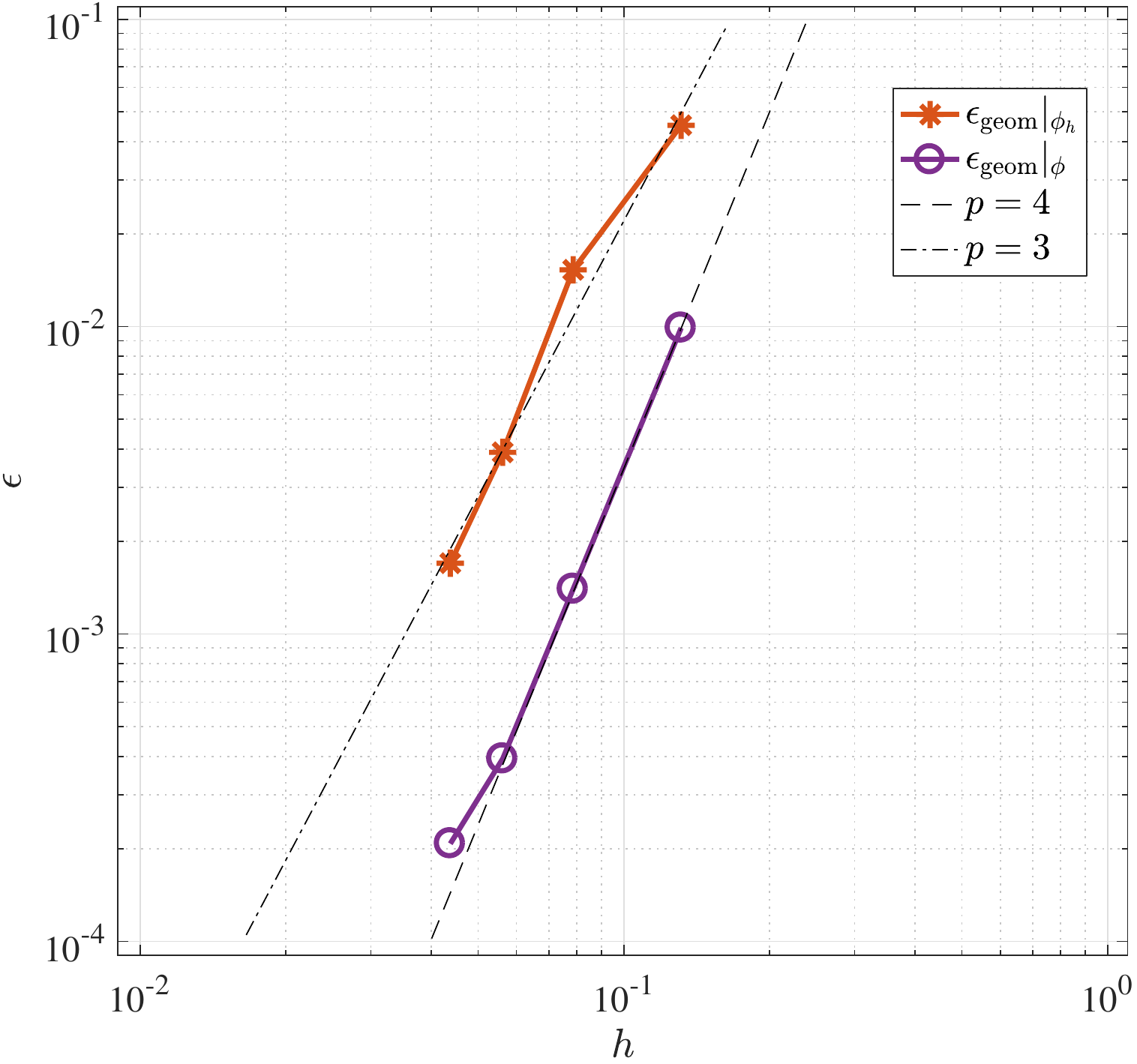}
\par\end{centering}
}$\quad$\subfloat[]{\begin{centering}
\includegraphics[width=0.49\textwidth]{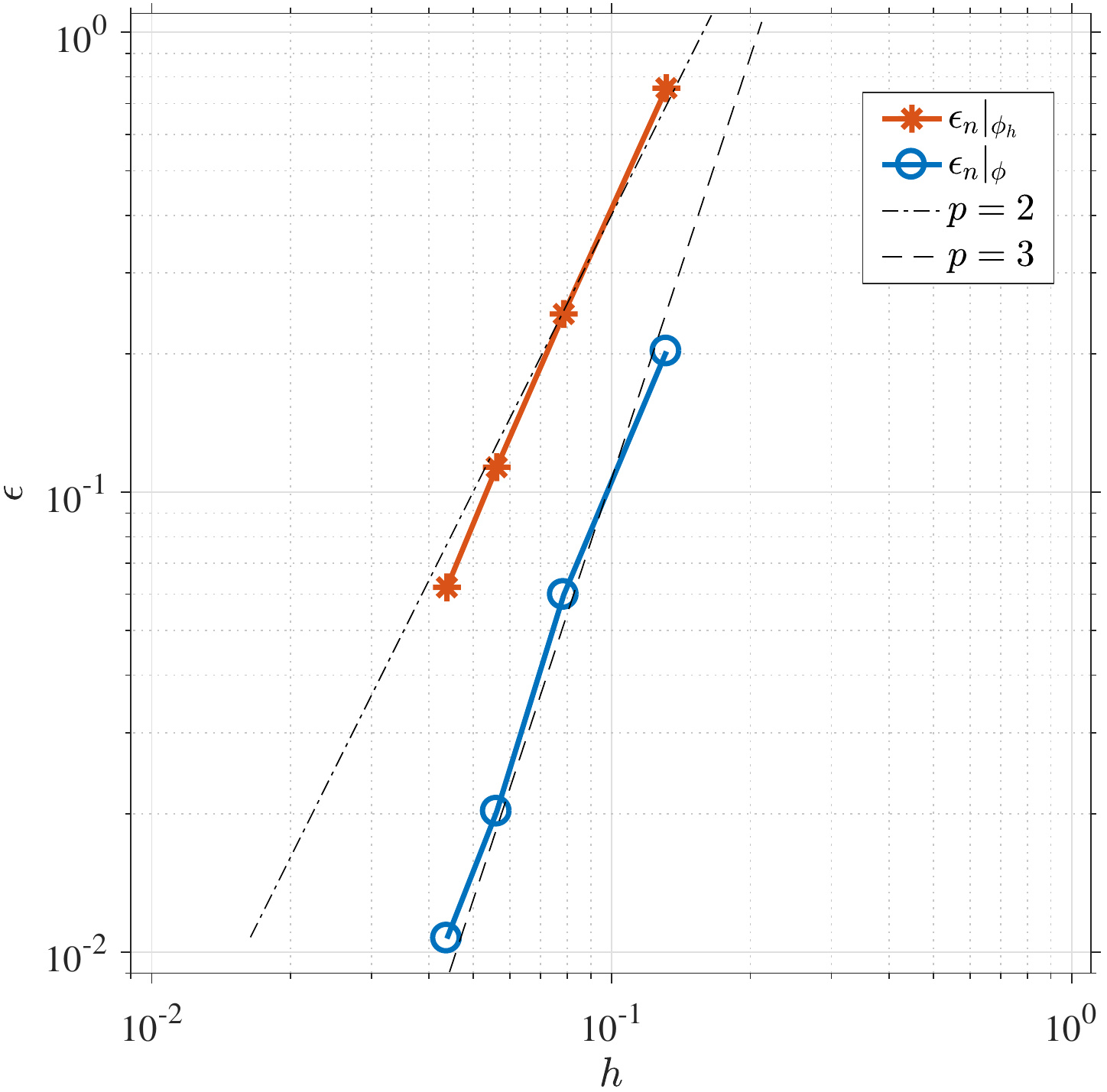}
\par\end{centering}
}\caption{Geometrical errors. (a) Distance error. (b) Normals error.\label{fig:Geometrical-errors}}
\end{figure}

In this section we numerically analyze the double approximation of
$\Gamma_{h}|_{\phi_{h}}$ compared to $\Gamma_{h}|_{\phi}$ by measuring
the distance error in $L_{2}$ -norm and computing the convergence
rates. The geometrical error with respect to the exact and discrete
distance function is given by

\[
\epsilon_{\mathrm{geom}}:=\|\phi(\bm{x}_{\Gamma_{h}})\|_{L_{2}(\Gamma_{h})},
\]
where $\bm{x}_{\Gamma_{h}}$ is the extracted surface using either
$\phi$ or $\phi_{h}$. Another interesting aspect with respect to
the tangential calculus approach is the normal errors introduced by
the discrete surface approximation. The normal error with respect
to the exact and discrete distance function is given by

\[
\epsilon_{n}:=\|\bm{n}_{e}-\bm{n}_{a}\|_{L_{2}(\Gamma_{h})},
\]
where $\bm{n}_{e}$ is the exact normal and $\bm{n}_{a}$ is the approximated
evaluated either with respect to $\phi$ or $\phi_{h}$. The results
can be found in Figure \ref{fig:Geometrical-errors} and Tables \ref{tab:Geometrical-errors-comparison.}
and \ref{tab:Normal-errors-comparison.}.

\begin{table}
\centering{}%
\begin{tabular}{|c|c|c|c|c|c|}
\hline 
$k$ & $h$ & $\epsilon_{\mathrm{geom}}|_{\phi}:=\|\phi(\bm{x}_{\Gamma_{h}}(\phi))\|_{L_{2}(\Gamma_{h})}$ & Rate & $\epsilon_{\mathrm{geom}}|_{\phi_{h}}:=\|\phi(\bm{x}_{\Gamma_{h}}(\phi_{h}))\|_{L_{2}(\Gamma_{h})}$ & Rate\tabularnewline
\hline 
\hline 
1 & 0.1314 & 0.0099 & - & 0.0452 & -\tabularnewline
\hline 
2 & 0.0786 & 0.0014 & 3.8065 & 0.0153 & 2.1080\tabularnewline
\hline 
3 & 0.0562 & 3.9275e-04 & 3.7890 & 0.0039 & 4.0747\tabularnewline
\hline 
4 & 0.0438 & 2.0799e-04 & 2.5500 & 0.0017 & 3.3309\tabularnewline
\hline 
\end{tabular}\caption{Geometrical errors comparison.\label{tab:Geometrical-errors-comparison.}}
\end{table}

\begin{table}
\centering{}%
\begin{tabular}{|c|c|c|c|c|c|}
\hline 
$k$ & $h$ & $\|\bm{n}_{e}-\bm{n}_{a}(\phi)\|_{L_{2}(\Gamma_{h})}$ & Rate & $\|\bm{n}_{e}-\bm{n}_{a}(\phi_{h})\|_{L_{2}(\Gamma_{h})}$ & Rate\tabularnewline
\hline 
\hline 
1 & 0.1314 & 0.2023 & - & 0.7562 & -\tabularnewline
\hline 
2 & 0.0786 & 0.0598 & 2.3717 & 0.2440 & 2.2012\tabularnewline
\hline 
3 & 0.0562 & 0.0202 & 3.2354 & 0.1133 & 2.2868\tabularnewline
\hline 
4 & 0.0438 & 0.0107 & 2.5491 & 0.0621 & 2.4121\tabularnewline
\hline 
\end{tabular}\caption{Normal errors comparison. \label{tab:Normal-errors-comparison.}}
\end{table}

\section{Concluding remarks\label{sec:Conclusion}}

In this paper we have introduced a finite element method for higher
order curved membranes using higher dimensional shape functions that
are restricted to the membrane surface. We have proposed a stabilization
for second order TraceFEM and show numerically that the solution is
stable and converges optimally. We have compared different parameterizations
and conclude that we get optimal convergence for the isoparametric
case $m_{B}=2$, $m_{\Gamma}=2$. We can observe that although no
solution stabilization is needed in the superparametric case of $m_{B}=1$,
$m_{\Gamma}=2$, with respect to mesh convergence, we still need stabilization
when interpolating the displacement field to the discrete surface,
cf. Figure \ref{fig:Displacement-fields}. The error difference between
the case of $m_{B}=1$, $m_{\Gamma}=1$ and $m_{B}=1$, $m_{\Gamma}=2$
is fairly small, cf. Figure \ref{fig:Stress-error-convergence}, and
since we still need a second order surface reconstruction for the
case of $m_{B}=1$, $m_{\Gamma}=2$, it seems natural to choose $m_{B}=2$,
$m_{\Gamma}=2$ instead. 

We have numerically shown the effects of different choices of the
stabilization parameters $\gamma_{1}$ and $\gamma_{2}$ and conclude
that the domain of optimal choices becomes bigger with smaller mesh
size.

The novelty of this work is the application of face stabilization
to second order TraceFEM for membrane problems. In future work we
will consider higher order TraceFEM using hexahedral elements.

\section*{Appendix}

\subsection*{Evaluation of basis functions in physical coordinates\label{subsec:Evaluation-of-basis}}

\begin{figure}
\centering{}\includegraphics[width=0.3\textwidth]{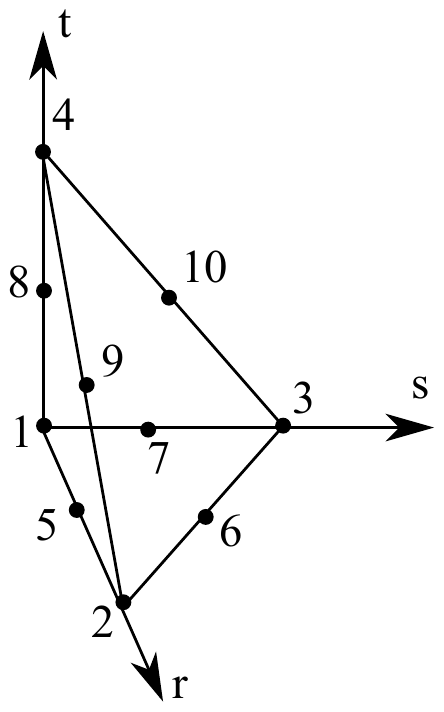}\caption{Tetrahedral node numbering.\label{fig:Tetrahedral-node-numbering.}}
\end{figure}

In order to construct $\varphi(\bm{x})$ on an affine second order
tetrahedron we define the geometric interpolation using the sub-parametric
mapping

\begin{equation}
\bm{x}=\sum_{i=1}^{4}\tilde{\varphi}_{i}\bm{x}_{i}\label{eq:subFi}
\end{equation}
where $\tilde{\varphi}_{i}$ are the basis function on the corner
nodes of a 10-noded tetrahedral element, with numbering according
to Figure \ref{fig:Tetrahedral-node-numbering.}, and $\bm{x}_{i}$
are the corresponding coordinates. We expand \ref{eq:subFi} and get

\[
(1-r-s-t)\bm{x}_{1}+r\bm{x}_{2}+s\bm{x}_{3}+t\bm{x}_{4}=\bm{x}
\]
which on matrix form is

\[
{\bf A\bm{r}+\bm{x}_{1}=\bm{x}}
\]
where 

\[
{\bf A=\begin{bmatrix}\bm{x}_{2}-\bm{x}_{1} & \bm{x}_{3}-\bm{x}_{1} & \bm{x}_{4}-\bm{x}_{1}\end{bmatrix}}
\]
with $\bm{x}_{i}=[x_{i},y_{i},z_{i}]^{\tr}$. We solve for $\bm{r}$
and get

\begin{equation}
\bm{r}(\bm{x})={\bf A^{-1}(\bm{x}-\bm{x}_{1}).}\label{eq:ReverseMap}
\end{equation}
Using the full basis function for the 10-noded tetrahedron $\varphi$
evaluated at $\bm{r}(\bm{x})$ we can write

\[
\varphi(\bm{x})=\varphi(\bm{r}(\bm{x}))
\]
and analogously 

\[
\nabla\varphi(\bm{x})=\nabla\varphi(\bm{r}(\bm{x})).
\]

Note that for every background element $K$, ${\bf A^{-1}}$ needs
only be computed once, which improves the performence of the root
finding method. 

\section*{Acknowledgement}
This research was supported by the Swedish Research Council Grant
No. 2011-4992.

\bibliographystyle{abbrv}
\bibliography{membraneP2ArXiv}

\end{document}